%% file: ABBS_Submitted.tex
\documentclass[11pt]{amsart}

\usepackage{a4wide,fullpage, amsfonts, amscd, amssymb}
\usepackage{epsf}
\usepackage{epsfig}
\usepackage{psfrag}
\usepackage{graphicx}
\newtheorem{theorem}{Theorem}[section]
\newtheorem{proposition}[theorem]{Proposition}
\newtheorem{lemma}[theorem]{Lemma}
\newtheorem{corollary}[theorem]{Corollary}
\newtheorem{definition}[theorem]{Definition}
\newtheorem{notation}[theorem]{Notation}
\newtheorem{theorem-definition}[theorem]{Theorem and Definition}

\theoremstyle{definition}

\newtheorem{remark}[theorem]{Remark}

\newcommand {\A}  {{\mathcal A}}

\newcommand {\Tt}{\widetilde{\mathcal {T}}}

\newcommand {\RR}  {{\mathbb R}}

\newcommand{\OO}{\mathfrak{O}}
\newcommand{\GN}{{\mathbb N}}
\newcommand{\GR}{{\mathbb R}}

\newcommand{\GZ}{{\mathbb Z}}
\newcommand{\GC}{{\mathbb C}}

\newcommand{\GK}{{\mathbb K}}
\newcommand{\GQ}{{\mathbb Q}}

\newcommand{\Int}{{\mathrm{Int}}}
\newcommand{\SC}{\mathcal{S}}
\newcommand{\Lb}{\mathcal{L}_\beta}
\newcommand{\La}[1]{{\mathcal{L}}_{\beta}^{(#1)}}

\newcommand{\Ta}[1]{\widetilde{\mathcal{T}}^{(#1)}}
\newcommand{\Sb}{\mathcal{S}_\beta}
\newcommand{\phb}{\phi_{\beta}}

\numberwithin{equation}{section}

\title[Pisot beta-numeration]{Boundary of central tiles
  associated with Pisot beta-numeration and purely periodic expansions} 

\dedicatory{Robert Tichy gewidmet, aus Anlass seines f\"unfzigsten
  Geburtstages.}

\author[S. Akiyama]{Shigeki Akiyama}
\address{Dept. of Mathematics, Faculty of Science Niigata University,
  Ikarashi-2, 8050, Niigata 950-2181, Japan}
\email{akiyama@math.sc.niigata-u.ac.jp}
\author[G. Barat]{Guy Barat}
\address{Institut f\"ur Mathematik A, T.U. Graz, Steyrergasse 30, 8010 Graz,
Austria}
\email{guy.barat@tugraz.at}
\author[V. Berthe\'e]{Val\'erie Berth\'e}
\address{LIRMM - CNRS UMR 5506-161 rue Ada, 34392 Montpellier Cedex 5, France} 
\email{berthe@lirmm.fr}
\author[A. Siegel]{Anne Siegel}
\address{IRISA - Campus de Beaulieu, 35042 Rennes Cedex, France}
\email{asiegel@irisa.fr}
\thanks{The second author was supported by the
Austrian Science Foundation FWF, project S9605, that is  part of the Austrian
National Research Network ``Analytic Combinatorics and Probabilistic
Number Theory''.}


\date{\today}

\begin{document}
\begin{abstract}This paper studies tilings related to  the
  $\beta$-transformation when $\beta$ is a Pisot number (that is not supposed
  to be a unit). Then it applies the obtained results to study the set of
  rational numbers  having a purely periodic  $\beta$-expansion.  Secial focus
  is given to some  quadratic examples.
\end{abstract}

\maketitle
%\tableofcontents

\section{Introduction}
Beta-numeration generalises usual binary and decimal numeration. Taking any
real number $\beta>1$, it consists in expanding numbers $x \in [0,1]$  as power
series in base $\beta^{-1}$ with digits in $\mathcal{D}=\{ 0, \dots,
\lceil\beta\rceil-1\}$.  As for $\beta\in\mathbb{N}$, the digits are obtained
with the so-called {\em greedy algorithm}: the $\beta$-transformation
$T_{\beta}: \ x \mapsto\beta x \pmod 1$ computes the digits $u_i=\lfloor \beta
T_{\beta}^{i-1}(x)\rfloor$, which yield the expansion $x =\sum_{i\geq 1} u_i
\beta^{-i}$. The sequence of digits is denoted  by $d_{\beta}(x)=(u_i)_{i\geq
  1}$.\\ 

The set of expansions $(u_i)_{i\geq 1}$ has been characterised by Parry  in
\cite{Parry60}) (see 
Theorem~\ref{parry} below). In case $\beta$ is a Pisot number,
Bertrand~\cite{BM77} and Schmidt~\cite{Schmidt80} proved independently that the
$\beta$-expansion $d_\beta(x)$ of a real number $x\in[0,1]$ is ultimately
periodic if and only if $x$ belongs to $\mathbb{Q}(\beta)\cap[0,1]$.  
A further natural question was to identify the set of numbers with purely
periodic expansions. For $\beta\in\GN$, it is known for a long time that
rational numbers $a/b$ with purely periodic $\beta$-expansion are exactly those
such that $b$ and $\beta$ are coprime, the length of the period being the order
of $\beta$ in $(\GZ/b\GZ)^*$.  Using an approximation and renormalisation
technique, Schmidt proved  in \cite{Schmidt80}
that when $\beta^2=n\beta+1$ and $n\in\GN^*$, then  all rational  numbers less
than 1 have a purely periodic $\beta$-expansion. This result was completed in~\cite{HamaImahashi},
to $\beta^2=n\beta-1$, $n\geq 3$, with 
respect to which no rational number has purely periodic
$\beta$-expansion. More generally, the latter result is satisfied by all
$\beta$'s admitting at least one positive real Galois conjugate in
$[0,1]$~\cite{Akiyama98}[Proposition~5].   Ito and Rao characterised the real
numbers having purely periodic $\beta$-expansion in terms of the associated
Rauzy fractal for any Pisot unit $\beta$ \cite{ItoRao04}, whereas the non-unit
case has been handled in \cite{BertheSiegel07}.  The length of the periodic
expansions with respect to quadratic Pisot units were investigated
in~\cite{QuRaoYang05}.\\  

Another natural question is to determinate the real numbers with finite
expansion. According to   \cite{FrougnySolomyak92}, we say that $\beta$
satisfies the  
finiteness property (F) if the
positive elements of $\mathbb{Z}[1/\beta]$ all have a finite $\beta$-expansion
(the converse is clear). A complete characterisation of $\beta$
satisfying the finiteness  property (F) is known when $\beta$ is a Pisot number
of degree 2 or 
3~\cite{Akiyama00}. It turns out that those numbers $\beta$ also play a role in
the question of purely periodic expansions. Indeed, if
$\beta$ is a unit Pisot number and satisfies the finiteness  property (F), then
there exists a neighbourhood of 0 in  $\mathbb{Q}_+$ whose elements all have
purely periodic $\beta$-expansion~\cite{Akiyama98}. This result is
quite unexpected since there is no reason {\em a priori} for obtaining only
purely periodic expansion around zero.\\  

The present paper investigates the case when $\beta$ is still a Pisot number,
but not necessarily a unit. We make use of the connection between pure
periodicity  and a compact self-similar representation of numbers having  no
fractional part in their $\beta$-expansion, similarly as in
\cite{ItoRao04,BertheSiegel07}. This representation  
is called the {\em central tile} associated with $\beta$ (\textit{Rauzy
  fractal}, or \textit{atomic surface} may also be encountered in the
literature, see e.g.  the survey \cite{BertheSiegel05}). For 
elements $x$ of the ring $\GZ[1/\beta]$, so-called $x$-tiles are introduced, so
that the central tile is a finite union of $x$-tiles up to translation. Those
$x$-tiles provide a covering of the space we are working in. We
first discuss the topological and metric properties of the central tile in the
flavor of~\cite{Akiyama02, Praggastis, Siegel03} and the relations between the
tiles.\\

In the unit case, the covering  by $x$-tiles is defined  in a Euclidean space
$\GK_\infty\simeq\GR^{r-1}\times\GC^s$, 
where $d=r+2s$ is the degree of the extension $[\GQ(\beta):\GQ]$ and $r$ the
number of real roots of beta's minimal polynomial. The space $\GK_\infty$ can be
interpreted as the product of all Archimedean completions of $\GQ(\beta)$
distinct from the usual one. It turns out that this is not enough in general:
in order to have suitable measure-preserving properties,
one has to take into  account the non-Archimedean completions associated with
the 
principal ideal $(\beta)$. Therefore, everything takes place in the product 
$\GK_\beta=\GK_\infty\times \GK_f$, where the latter is a finite product of
local fields. In the framework of substitutions, 
this approach has been already used in~\cite{Siegel03}, and  was inspired by
\cite{rau2}.  See also \cite{Sing06}. Completions and (complete) tiles are
introduced in Section~\ref{sec:complete}. We discuss why taking  into  account
non-Archimedean completions 
is  suitable from a tiling point of view: when the finiteness property (F)
holds, we prove that the $x$-tiles are disjoint 
if the non-Archimedean completions are considered, which was not the case when
only taking into account Archimedean completions.  
Our principal result in that context is Theorem~\ref{thm:disjointness}.

Our main goal is the study of the set of rational numbers having purely periodic
beta-expansion, for which we introduce the following notation.

\begin{notation}\label{def:periodiques} $\Pi_\beta$ denotes the set of those real numbers
  $x\in[0,1)$ having purely periodic
  beta-expansion. We also note $\Pi_\beta^{(r)}=\Pi_\beta\cap\GQ$.
\end{notation}
The study of those sets begins in
Section~\ref{sec:periodicity}. After having recalled the characterisation of 
purely periodic expansions in terms of the complete
tiles due to~\cite{BertheSiegel07} (see~\cite{ItoRao04} for the unit case), we
apply it to gain results on the periodic expansions of the rational integers. 

\begin{theorem}
Let $\beta$ be a Pisot number that satisfies the property $(F)$. Then there
exist $\varepsilon$ and $D$ such that  for every $x=\frac{p}{q} \in {\mathbb
  Q}\cap [0,1)$, if $x \leq \varepsilon$, $\gcd(N(\beta),q)=1$ and 
$N(\beta)^D$ divides $p$, then $x$ has a purely periodic expansion in base
$\beta$. 
\end{theorem}
\begin{definition}[Function gamma] \label{def:gamma}
The function $\gamma $ is defined  on the set of Pisot numbers  and takes its
values in $[0,1]$. Let $\beta$ be a Pisot number. Let $N(\beta)$ denote the
norm of $\beta$. Then, $\gamma(\beta)$ is defined as  
\[
\gamma(\beta) = \sup\left\{ v \in [0,1]; \  \forall x=\frac{p}{q}
  \in {\mathbb Q}\cap]0,v] \mbox{ with}  
  \gcd(q,N(\beta))=1, \mbox{ then } x\in\Pi_\beta^{(r)}
\right\}. 
\] 
\end{definition}
The reasons of the condition $\gcd(q,N(\beta))=1$ will be given in
Lemma~\ref{lem:pp}. We also use the central tile and its tiling properties to
obtain in Section~\ref{sec:examples} an explicit  
computation of the quantity $\gamma(\beta)$ for two quadratic Pisot numbers,
that is, 
\begin{theorem} 
 $\gamma(2+ \sqrt 7)=0$ and   $\gamma(5+2\sqrt 7)=(7-\sqrt 7)/12$. 
\end{theorem}

The second example shows that the behaviour of $\gamma(\beta)$ in the non-unit
case is slightly different from its behaviour in  the unit case. \\ 

This paper is organised as follows. Section~\ref{sec:vocabulaire} recalls
terminology and results necessary to 
state and to prove the results, including Euclidean tiles and the unit
case. Section~\ref{sec:complete} goes beyond the unit case and extends the
previous concepts including non-Archimedean components. This section starts
with a short compendium on what we need from algebraic number
theory. Section~\ref{sec:periodicity} studies purely periodic expansions and
the Section~\ref{sec:examples} is devoted to examples in quadratic fields.\\

Since we work with Pisot numbers and in order
to avoid the introduction of plethoric vocabulary, we will always assume in
this section that $\beta$ is a Pisot number, even if the result is more
general. The reader interested in generalities concerning beta-numeration could
have a look to \cite{Blanchard, BertheSiegel05, Baratetal}.

\section{Beta-numeration, automata, and tiles}\label{sec:vocabulaire} 

\subsection{Beta-numeration}

We assume  that $\beta$ be a Pisot number.
Since $1\in\mathbb{Q}$,  $d_{\beta}(1)$ is ultimately periodic by
\cite{BM77,Schmidt80} and we have the following (see e.g. \cite{Parry60,
  Blanchard, FrougnyTemuco, FrougnyLothaire}):

\begin{theorem-definition}\label{parry} 
Let $\beta$ be a Pisot number. Let $\mathcal{D}=\{0,1,\cdots, \lceil
\beta\rceil -1\}$. Let $d^*_{\beta}(1)=d_{\beta}(1)$ if $d_{\beta}(1)$ 
is infinite and $d^*_{\beta}(1)=(t_1 \dots t_{n-1} t_n)^{\infty}$,
if $d_{\beta}(1)=t_1\dots t_{n-1}(t_n+1)0^\infty$, with $t_i \in  \mathcal{D}$
for all $i$. Then 
the set  of $\beta$-expansions of real numbers in $[0,1)$
is exactly the set  of sequences $(u_i)_{i\geq 1}$ in ${\mathcal D}^{\GN}$
that satisfy the so-called {\em admissibility condition}
\begin{equation}\label{eq:real}
\forall k \geq 1, \ (u_i)_{i\geq k} <_{\mbox{lex}} d^*_{\beta}(1).
\end{equation}
A finite string $w$ is said to be admissible if the sequence $w\cdot 0^\infty$
satisfies the condition~(\ref{eq:real}), where $A\cdot B$ denotes the
concatenation of the words $A$ and $B$. The set of admissible strings is
denoted  by $\Lb$; the set of admissible sequences by $\Lb^{\infty}$. The map
$x\mapsto d_\beta(x)$ realises an increasing bijection from $[0,1)$ onto $\Lb$,
endowed with the lexicographical order.\\ 
\end{theorem-definition}

\begin{notation}\label{not:rn}
From now on, $\beta$ will be a Pisot number of degree $d$, with $$d_\beta^*(1)=
t_1\cdots t_m(t_{m+1}\cdots t_{n})^\infty,$$ that is, $n$ is the 
sum of the lengths of the preperiod and of the period; in particular, $m=0$ if
and only if  $d_\beta^*(1)$ is purely periodic .
\end{notation}
  The Pisot number $\beta$
  is said to be a  {\em simple Parry number} if
  $d_{\beta}(1)$ is  finite, it is said to be a  {\em non-simple Parry
    number}, otherwise. 
  One has $m=0$  if and only  if $\beta$ is a simple Parry number: indeed,
  $d_{\beta}(1)$ is never purely periodic 
 according to Remark  7.2.5 in \cite{FrougnyLothaire}). We will denote by
  ${\mathcal A}$ the alphabet $\{1, \dots, n\}$.\\

\paragraph{\bf Expansion of the non-negative real numbers}
The $\beta$-expansion of any $x \in {\mathbb R}^+$ is deduced by rescaling from
the expansion of $\beta^{-p}x$, where 
$p$ is the smallest integer such that $\beta^{-p}x \in [0,1)$: 
\begin{equation}\label{eq:expansion}
\forall \, x \in {\mathbb R}^+, \, x=\underbrace{w_{p}\beta^p + \dots +
  w_0}_{\mbox{integer part}} + \underbrace{u_1 \beta^{-1} + \cdots + u_i^i
  \beta^{-i} + \cdots}_{\mbox{fractional part}}, \quad w_p\cdots w_0u_1 \cdots
u_i \cdots \;\textrm{satisfies}~(\ref{eq:real}).
\end{equation}

In this case, we call $[x]_\beta=w_{p}\beta^p + \dots + w_0 $ the {\em integer
  part} of $x$  and $\{x\}_\beta=u_1 \beta^{-1} + \cdots + u_i^i \beta^{-i} +
\cdots$ the {\em fractional part} of $x$. We extend the notation $d_\beta$ and
write $d_\beta(x)=w_p\cdots w_0. u_1 \cdots u_i \cdots$. 

\medskip

\paragraph{\bf Integers in base $\beta$}
We define the set of integers in base $\beta$ as the set of positive real
numbers with no fractional part:
\begin{align}
\Int(\beta)&= \left\{ w_p \beta^p + \dots + w_0  ; \,   w_p \dots w_0
  \in\Lb\right\}\\ 
&=\left\{ [x]_\beta;\, x\in \mathbb{R}_+\right\}\subset{\mathbb Z}[\beta].\notag
\end{align}

The set $\Int(\beta)$ builds a discrete subset of $\mathbb{R}_+$. It has some
regularity: two consecutive 
points in $\Int(\beta)$  differ  by a finite number of values, namely, the
positive numbers $T_{\beta}^{a-1}(1)$, $a \in  \{1,\cdots,
n\}$ (see \cite{Thurston,Akiyama07}). It can even be shown that it is a Meyer
set~\cite{BFGK}.

\subsection{Admissibility graph}\label{sec:ag}

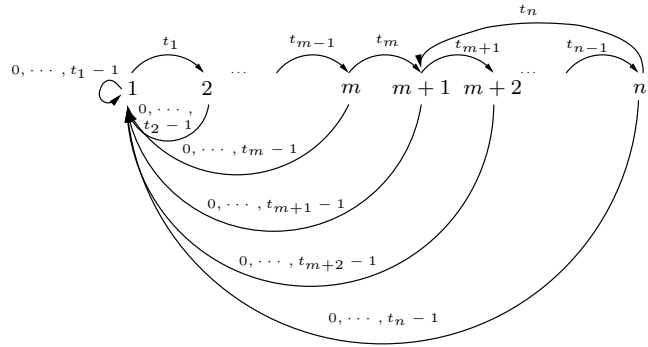
\begin{figure}[ht]
\begin{center}
\input{graph2.pstex_t}
\caption{
 The graph describes admissible sequences for the $\beta$-shift. The
 number $n$ of nodes is given by  the sum of the preperiod and the period of
$d^*_{\beta}(1)=t_1\cdots t_m(t_{m+1}\cdots t_{n})^\infty$. 
 From each node $a$ to the node $1$, there are $t_a$ edges labelled by $0,
\dots, t_a-1$. From each node ${a}$ to the node $a+1$,  there is one  edge
labelled by  $t_{a}$. Let $m$ denote the length of the preperiod of
$d^*_{\beta}(1)$ (it can  possibly be zero). From the node $n$ to the node
$m+1$ there is an edge labelled by $t_n$.
}\label{Fig:admissibility}
\end{center} 
\end{figure}

The set of admissible sequences described by~(\ref{eq:real}) is the set of
infinite labellings of an explicit finite graph with nodes in  
${\mathcal A}=\{1, \dots, n\}$  and edges $b\xrightarrow{\varepsilon}a$, with
$a,b \in \mathcal{A}$  labelled by  digits
$\varepsilon \in  \mathcal{D}=\{0,1,\cdots, \lceil \beta\rceil -1\}$. This
so-called {\em admissibility graph} is 
depicted  in Figure~\ref{Fig:admissibility}.

For $a\in\A$, define $\mathcal{L}_{\beta} ^{(a)}$  as  the set
of admissible strings  $w$ (see Definition  \ref{parry})   that the graph of
admissibility conducts from the 
initial node 1 to the node $a$. In other words, for $a\neq 1$,
$\mathcal{L}_{\beta} ^{(a)}$ is 
the set of admissible strings having $t_1\cdots t_{a-1}$ as a suffix. Clearly,
according to the form of the admissibility  graph, one has 
$\Lb=\bigcup_{a\in\mathcal A}\mathcal{L}_{\beta}^{(a)}$.\\

Denote by $S$ the shift operator on the  set of sequences  in the set of digits
$\{0,\ldots, \lceil\beta\rceil-1\}^\GN=\mathcal{D}^\GN$. The beta-expansion of
$T_\beta^k(1)$ is 
$d_\beta(T_\beta^k(1))=S^k(d_\beta(1))$. By increasingness of the map
$d_\beta$, it follows that for any $x\in[0,1)$:
\begin{align}\label{T(a)}
t_1t_2\cdots t_{a-1} d_\beta(x)\in\mathcal{L}_{\beta}
^{\infty}&\Longleftrightarrow 
d_\beta(x)<_{\rm lex}S^{a-1}(d_\beta^*(1))\\
&\Longleftrightarrow x\in[0,T^{a-1}(1)).\notag
\end{align}
Notice that if $\beta$ is a simple Parry number (that is, if $m=0$) and
$k\in\GN$,  then the sequence $S^k(d_\beta^*(1))$ is not admissible.

\subsection{Central tiles}

The {\em central tile} associated with a Pisot number is a compact geometric
representation of the set $\Int(\beta)$ of integers in base $\beta$. It is
defined as follows. 

\medskip

\paragraph{\bf Galois conjugates of $\beta$ and euclidean completions}
Let  $\beta_2$, $\dots$, $\beta_r$  be the 
real conjugates of $\beta=\beta_1$; they all have modulus strictly smaller than
$1$, since $\beta$ is a Pisot number. Let ${\beta}_{r+1}$,
$\overline{\beta_{r+1}}$, 
$\dots$, $\beta_{r+s}$,  $\overline{\beta_{r+s}}$ stand for its complex
conjugates. For $2 \leq i \leq  r$, let ${\mathbb K}_{\beta_i}$ be equal
to ${\mathbb R}$, and for $r+1\leq i \leq
r+s$, let ${\mathbb K}_{\beta_i}$ be equal
to  ${\mathbb C}$. The fields ${\mathbb R}$  and ${\mathbb C}$  are endowed
with the normalised absolute value $\vert x\vert_{\mathbb{K}_{\beta_i}}=\vert
x\vert$ 
if $\mathbb{K}_{\beta_i}=\mathbb{R}$ and $\vert x\vert_{\mathbb{K}_{\beta_i}}=\vert
x\vert^2$ if $\mathbb{K}_{\beta_i}=\mathbb{C}$. Those absolute values induce
the usual topologies on $\mathbb{R}$ (resp. $\mathbb{C}$). 
For any $i=2$ to $r+s$, the  $\GQ$-homomorphism defined on
$\GQ(\beta)$ by $\tau_i(\beta)=\beta_i$ realises a $\GQ$-isomorphism between
$\GQ(\beta)$ and $\GK_{\beta_i}=\GQ(\beta_i)\hookrightarrow\mathbb{R},
\mathbb{C}$.\\

\paragraph{\bf Euclidean  $\beta$-representation space}
We obtain a Euclidean representation $\mathbb{Q}$-vector space
$\mathbb{K}_\infty$ by gathering the fields  ${\mathbb K}_{\beta_i}$:  
$${\mathbb K}_{\infty}={\mathbb K}_{\beta_2}\times \dots\times{\mathbb
  K}_{\beta_r}\times {\mathbb K}_{\beta_{r+1}}\times{\mathbb
  K}_{\beta_{r+2}}\times\cdots\times {\mathbb K }_{\beta_{r+s}} \simeq
\GR^{r-1} \times\GC^s.$$ 

We denote by $\Vert\cdot\Vert_{\infty}$ the maximum norm on $\GK_{\infty}$. We
have a natural embedding 
\begin{equation*}
\begin{array}{rcl}
\phi_\infty\colon\mathbb{Q}(\beta)&\longrightarrow&\mathbb{K}_\infty\\
x&\longmapsto&(\tau_i(x))_{2\leq i\leq r+s}
\end{array}
\end{equation*}

\medskip
 
\paragraph{\bf Euclidean central tile}
We are now able to define the central tiles and its associated subtiles:

\begin{definition}[Central tile]
Let $\beta$ be a Pisot number with degree $d$. The {\em Euclidean central tile}
of $\beta$ is the representation of the set of integers in base $\beta$: 
\[
{\mathcal T} =  \overline{\phi_{\infty}(\Int(\beta))} \subset
\overline{\GQ(\beta_2)\times\cdots\times\GQ(\beta_{r+s})}\subset\GK_{\infty}.
\]
\end{definition}
Since the roots $\beta_i$ have modulus smaller than one, $\mathcal T$ is a
compact subset of $\GK_{\infty}.$

\subsection{Property (F) and tilings}
More generally, to each $x \in {\mathbb Z}[1/\beta]\cap[0,1)$  we associate a
geometric representation of points that admit $w$ as fractional part. 

\begin{definition}[$x$-tile]
Let $x \in {\mathbb Z}[1/\beta]\cap[0,1)$.   The tile associated with $x$  is 
$${\mathcal T}(x) =  \overline{ \phi_{\infty}(\{y \in {\mathbb R}^+; \,
  \{y\}_\beta=x \}) } \subset \phi_{\infty}(x) + {\mathcal T}.$$  
\end{definition} 

It is proved in \cite{Akiyama02} that the tiles
${\mathcal T}(x)$ provide a covering of ${\mathbb K}_{\infty}$, i.e.,

\begin{equation} \label{eq:covunit}
\GK_{\infty}= \bigcup_{x \in\GZ [1/\beta]\cap[0,1)}\mathcal{T}(x).
\end{equation}

Since we know that the tiles ${\mathcal T}(x)$ cover the space ${\mathbb
  K}_{\infty}$, a natural question is whether this covering is a tiling (up to
sets of zero measure).  
  \begin{definition}[Exclusive points] \label{def:exclu}We say
that a point $z \in {\mathbb K}_{\infty}$ is {\em exclusive} in the tile
${\mathcal T}(x)$ if $z$  is  contained in  no other tile ${\mathcal
  T}(x^\prime)$ with $x^\prime \in {\mathbb Z}[1/\beta]\cap[0,1)$,   and
$x'\neq x$.  
  \end{definition}
\begin {definition}[Finiteness property]\label{def:F}
 The Pisot number $\beta$ satisfies the {\em  finiteness property (F)} if and
 only if every $x \in {\mathbb Z}[1/\beta]\cap[0,1)$ has a finite
 $\beta$-expansion.
\end{definition}
If the finiteness property is satisfied, a sufficient  tiling condition  is
known when $\beta$ is a unit.   

\begin{theorem}[Tiling property]\label{thm:Pavageunit}
Let $\beta$ be a unit Pisot number. The number  $\beta$ satisfies the
finiteness property (F) if 
and only if $0$ is an exclusive inner point of the central tile of $\beta$. 
In this latter case,  every  tile  $\mathcal{T}(x)$,  for $x\in\GZ
[1/\beta]\cap[0,1]$ 
has a non-empty interior, and all its inner points  are exclusive. In other
words,  the tiles  $\mathcal{T}(x)$ provide a tiling of $\GK_{\infty}$.
\end{theorem}
\begin{proof}  The proof  is
  done in~\cite{Akiyama02}.  In~\cite{SiegelThuswaldner}, this property is
  restated in a discrete geometry framework. 
\end{proof}

\subsection{Purely periodic points}

In~\cite{ItoRao04}, Ito and Rao establish a relation between the central tile
and purely periodic $\beta$-expansions. For that purpose, a geometric
realisation of the natural extension of the beta-transformation is built using
the central tile.  More precisely, the central tile represents by construction
(up to closure) the  strings $w_n \dots w_0$ that can be read in the
admissibility graph shown in Figure~\ref{Fig:admissibility}. We gather strings
$w_m \dots w_0$ depending on the nodes of the graph to which the string
$w_m\dots w_0$ arrives. 

\begin{definition}[Central subtiles]
Let $a\in \mathcal{A}=\{1,\cdots,n\}$.
The {\em central $a$-subtile} is defined as 
\[
{\mathcal T}^{(a)} =\overline{
  \phi_{\infty}\left(\left\{x\in\Int(\beta);\,
  d_\beta(x)\in\La{a}\right\}\right)}.
\]
\end{definition}

\begin{theorem}[\cite{ItoRao04}]\label{thm:periodicunit}
Let $\beta$ be a Pisot unit.  We recall that  $\mathcal{A}=\{1,\cdots,n\}$.
 Let $x \in {\mathbb Q}(\beta)\cap[0,1)$. The 
$\beta$-expansion of $x$ is purely periodic if and only if
$$(-\phi_{\infty}(x),x) \in \bigcup_{a \in \A}{\mathcal
  T}^{(a)} \times[0, T_{\beta}^{a-1}(1)).$$
\end{theorem}

As soon as 0 is an inner point of the central tile, we deduce that small
rational numbers  have a purely periodic expansion.  

\begin{corollary}[\cite{Akiyama98}]\label{gamabetaunit}
Let $\beta$ be a Pisot unit. If $\beta$ satisfies  the finiteness property (F),
then there 
exists a constant $c>0$ such that every $x \in {\mathbb Q}\cap[0,c)$  has a
purely periodic expansion in base $\beta$. 
\end{corollary}
\begin{proof}
Since $0$ is an inner point of $\mathcal{T}$ and $\A$ is finite, there
exists $c>0$ such that $0<c\leq \min\{T_\beta^{(a-1)}(1);\, a\in\A\}$
and $B_\infty(0,c)\subset\mathcal{T}$. For $x\in[0,c)$, we have
$\phi_\infty(x)=(x,x,\ldots,x)$ and 
\[
(-\phi_\infty(x),x)\in\mathcal{T}\times[0,c)
\subset\bigcup_{a\in\A}\left(\mathcal{T}^{(a)}
  \times\left\lbrack 0,T_\beta^{(a-1)}(1)\right)\right).
\] 
Then the periodicity follows from Theorem~\ref{thm:periodicunit}.
\end{proof}

This result was first proved directly by Akiyama~\cite{Akiyama98}.
Recall that $\gamma(\beta)$ is the supremum of such $c$'s according to
Definition \ref{def:gamma}.  As soon 
as one of the  conjugates of $\beta$ is positive, then $\gamma(\beta)=0$.
The  quadractic unit case  is completely understood: in this  case  
Ito and Rao proved that $\gamma(\beta)$ equals 0 or
1~(\cite{ItoRao04}). Examples of computations of $\gamma(\beta)$  for higher
degrees are also  performed by Akiyama in the unit case in~\cite{Akiyama98}. 

\medskip

\paragraph{\bf Algebraic natural extension} By abuse of language, one may say that 
 Theorem \ref{thm:periodicunit} implies that $\bigcup_{a \in \A}
(\mathcal{T}^{(a)}\times\big[0,T_{\beta}^{a-1}(1)\big)$ is a fundamental domain
for  an   {\em algebraic} realisation of the  natural extension of the
$\beta$-transformation $T_{\beta}$, though it does not satisfy  Rohklin's
minimality condition  for natural extensions 
(see \cite{Rohlin}  and  also  \cite{CFS}). We wish to explain shortly this
reason in the sequel. 

In \cite{Solomyak}, Dajani~\textit{et al.} provide  an explicit  construction
of  the natural extension of  the $\beta$-transformation for any 
$\beta>1$ in dimension three,  the third dimension  being given by the
height in  a  stacking structure. This construction is minimal in the above
sense. As a by-product, one can retrieve the invariant measure of the system as
an induced measure. 
However, this natural extention  provides  no information
on the purely  periodic orbits under the action of the   $\beta$-transformation
$T_{\beta}$. The essential reason is that the  
geometric realisation map which plays the role of  our  $\phi_{\infty}$ is  not
an additive homomorphism. 
And therefore, this embedding destroys the algebraic structure of the 
$\beta$-transformation. 
Our construction, which  was originated by Thurston in the Pisot unit case
\cite{Thurston},   only works for  restricted cases but it has 
the advantage  that
we can use the conjugate maps which are additive homomorphisms. 
This is the clue  used by  Ito and Rao  in~\cite{ItoRao04} for the
description of purely  periodic orbits. Summing up, we need a more geometric
natural extension than that of Rohklin 
to answer number theoretical questions like  periodicity issues.

Let us note that  in the non-unit case,  measure-preserving properties are no
more satisfied  by  the embedding  $\phi_{\infty}$. Indeed, it is clear that
$T_{\beta}$ is an expanding map with ratio $\beta$. By  involving only
Archimedean embeddings as in the unit case,  we will  only take into account
$\phi_{\infty}$ which  is a contracting  map with ratio $N(\beta)/\beta$ and we
won't be able to  get a  measure-preserving  natural extension.   This is the essential reason why we  introduce now non-Archimedean embeddings.

\section{Complete tilings}\label{sec:complete}

Thanks to the  non-Archimedean part, we will show that we obtain a map $\phb$
which is a contracting map with ratio $1/\beta$. Let us recall that $T_{\beta}$
is an expanding map with ratio $\beta$.  We thus  will recover a  realisation
of the natural extension via a  measure-preserving map. 
Moreover, the   extended map acting on the fundamental domain of  the
natural extension  will be almost one-to-one (being a  kind of variant of
Baker's transform). 
Therefore we have good chances to have a one-to-one map on the 
{\em lattice points}   for this algebraic
natural extension. Considering that a bijection
on a finite set yields  purely periodic expansions, we will obtain a
description of purely periodic elements of this system. This heuristics  will
be in fact realised in Proposition~\ref{prop:covering} and Theorem~\ref{thm:disjointness} below.

\subsection{Algebraic framework}\label{subsec:bourbaki}

In order to extend the results above to the case
where $\beta$ is not a  unit, we follow the idea of~\cite{Siegel03} and embed
the  central tile in a larger space including local components. To avoid
confusion, the central tile ${\mathcal T}\subset \RR^{r-1}\times\GC^s$ will be
called the {\em  Euclidean central tile}.  The large tile will be called {\em
  complete tile} and  denoted as  $\widetilde{\mathcal T}$. \\ 

Let us briefly recall some facts and set notation. The results can be found for
instance in the first two chapters of~\cite{CasselsFroehlich}. Let
$\mathfrak{O}$ be the ring of 
integers of the field $\mathbb{Q}(\beta)$. If $\mathfrak{P}$ is a prime ideal
in $\OO$ such that $\mathfrak{P}\cap\GZ=p\GZ$, with relative degree
$f(\mathfrak{P})=[\OO/\mathfrak{P}:\mathbb{Z}/p\GZ]$ and ramification index
$e(\mathfrak{P})$, then $\GK_{\mathfrak{P}}$ stands for the completion of
${\mathbb Q}(\beta)$ with respect to the ${\mathfrak P}$-adic topology. It is
an extension of $\mathbb{Q}_p$ of degree $e(\mathfrak{P})f(\mathfrak{P})$. The
corresponding normalised absolute value is given by $\vert
x\vert_{\mathfrak{P}}=\left\vert
  N_{\GK_\mathfrak{P}/\GQ_p}(y)\right\vert_p^{1/e(\mathfrak{P})f(\mathfrak{P})} =
p^{-f(\mathfrak{P})v_\mathfrak{P}(y)}$.  
We denote $\OO_\mathfrak{P}$ its ring of integers and
$\mathfrak{p}_{\mathfrak{P}}$ its maximal ideal; then 
\begin{gather*}
\OO_\mathfrak{P}=\{y\in\GK_\mathfrak{P}\, ;\, v_\mathfrak{P}(y)\geq
0\}=\{y\in\GK_\mathfrak{P}\, ;\, \vert y\vert_\mathfrak{P}\leq 1\}.\\
\mathfrak{p}_{\mathfrak{P}}=\{y\in\GK_\mathfrak{P}\, ;\, v_\mathfrak{P}(y)\geq
1\}=\{y\in\GK_\mathfrak{P}\, ;\, \vert y\vert_\mathfrak{P}< 1\}.
\end{gather*}
The normalised Haar measure on $\GK_{\mathfrak{P}}$ is
$\mu_{\mathfrak{P}}(a+\mathfrak{p}_{\mathfrak{P}}^m)= p^{mf(\mathfrak{P})}$. In
particular: $\mu_{\mathfrak{P}}(\OO_\mathfrak{P})=1$.

\begin{lemma}\label{lem:approximation} Let $\mathcal{V}$ be the set of places
in $\mathbb{Q}(\beta)$. For any place $v\in\mathcal{V}$, the
associated normalised absolute value is denoted $\vert \cdot\vert_v$. If $v$ is
Archimedean, we make the usual convention $\OO_v=\GK_v$.\\
$(1)$ Let $\SC\subset\mathcal V$ a finite set of places. Let $(a_v)_{v\in
  \SC}\in\prod_{v\in \SC}\GK_v$. Then, for any $\varepsilon>0$, there exists
$x\in\GK$ such that $\vert x-a_v\vert_v\leq\varepsilon$ for all $v\in \SC$.\\
$(2)$ Let $\SC\subset\mathcal V$ a finite set of places and
$v_0\in\mathcal{V}\setminus \SC$. Let $(a_v)_{v\in
 \SC}\in\prod_{v\in \SC}\GK_v$. Then, for any $\varepsilon>0$, there exists
$x\in\GK$ such that $\vert x-a_v\vert_v\leq\varepsilon$ for all $v\in \SC$ and
$v\in\OO_s$ for all $v\not\in \SC\cup\{v_0\}$. Furthermore, if $v_0$ is an
Archimedean place and  $(a_v)_{v\in \SC}\in\prod_{v\in \SC}\OO_v$, then $x\in\OO$.
\end{lemma}
\begin{proof} (1) (resp.  the first part of (2)) are widely known as the weak
  (resp. strong) approximation theorems. Concerning the last sentence, let
  $x\in\GQ(\beta)$  given by (2). By assumption, $x\in\OO_v$ for all $v$,
  therefore $x\in\OO$, since $\OO$ is the intersection of the local rings
  $\OO_v$, where $v$ runs along the non-Archimedean places. 
\end{proof}

\subsection{Complete representation space}\label{subsec:crs}
\begin{notation}
Let $\mathfrak{P}_1$, $\dots$,  $\mathfrak{P}_{\nu}$ be  the prime ideals in
the ring of integers $\OO$ that contain $\beta$, that is, 
\[
(\beta)=\beta \OO=\prod_{i=1}^\nu {\mathfrak{P_i}}^{n_i}.
\]
For $x\in\GQ(\beta)$, $N(x)$ shortly denotes the norm $N_{\GQ(\beta)/\GQ}(x)$.
We have $N(\beta \OO)=\vert N(\beta)\vert$; the prime numbers 
$p$ arising from $\mathfrak{P}_i\cap \GZ=p\GZ$ are the prime factors of
$N(\beta)$. Let $\Sb$ be the set
containing the Archimedean places corresponding to $\beta_i$, $2\leq i\leq r+s$
and the $\nu$ non-Archimedean places corresponding to the $\mathfrak{P}_i$. 
\end{notation}

The {\em complete representation space} $\GK_\beta$ is obtained by adjoining to
the Euclidean representation the product of local fields
$\GK_f=\prod_{i=1}^\nu\GK_{\mathfrak{P}_i}$, that is
$\GK_\beta=\GK_\infty\times\GK_f=\prod_{v\in\Sb}\GK_v$. The field $\GQ(\beta)$
naturally embeds in $\GK_\beta$: 
\begin{equation*}
\begin{array}{rcl}
\phb\colon\mathbb{Q}(\beta)&\longrightarrow&\displaystyle{
  \mathbb{K}_\infty\times\prod_{i=1}^\nu\GK_{\mathfrak{P}_i}}\\   
x&\longmapsto&(\phi_\infty(x), x,\ldots, x)
\end{array}
\end{equation*}
The complete representation space is endowed with the product topology, and
with coordinatewise addition and multiplication. This makes it a locally
compact abelian ring. Then the approximation theorems yield the following:

\begin{lemma}\label{lem:density} With the previous notation, we have that
  $\phb(\GQ(\beta))$ is dense in $\GK_\beta$, and that $\phb(\OO)$
  is dense in $\prod_{v\in\Sb}\mathfrak{O}_v$.
\end{lemma}
\begin{proof}The first assertion follows from the first part of
  Lemma~\ref{lem:approximation} with $\SC=\Sb$. The second assertion
  follows from its second part with $\SC=\Sb$ and $v_0$ being the
  Archimedean valuation corresponding to the trivial embedding
  $\tau(\beta)=\beta$. 
\end{proof}

The normalised Haar measure $\mu_\beta$ of the additive group $(\GK_\beta,+)$
is the  product measure of the normalised Haar measures on the complete fields
$\GK_{\beta_i}$ (Lebesgue measure) and $\GK_{\mathfrak{P}_i}$ (Haar measure
$\mu_{\mathfrak{P}_i}$). By a standard measure-theoretical argument, if
$\alpha\in\GQ(\beta)$ and if $B$ is a borelian subset of $\GK_\beta$, then 
\begin{equation}\label{equ:homothetie}
\mu_\beta(\alpha\cdot B)
=\mu_\beta(B)\prod_{v\in\Sb}\vert\alpha\vert_v.  
\end{equation}
Consequently, if $\alpha\in\GQ(\beta)$ is a $\Sb$-unit (that is, if
$\vert\alpha\vert_v=1$ for all $v\not\in\Sb$), then
$\mu_\beta(\alpha\cdot B)=\vert\alpha\vert^{-1}\mu_\beta(B)$ by the product
formula ($\vert\cdot\vert$ is there the usual real absolute value). This holds
in particular for $\alpha=\beta$.\\

At last, we also denote by $\Vert\cdot\Vert$ the maximum norm on
${\mathbb K}_{\beta}$, that is $\displaystyle{\Vert
x\Vert=\max_{v\in\Sb}\vert x\vert_v}$. The following
finiteness remark will be used several times. 

\begin{lemma}\label{lem:finite}If $B\subset\GK_\beta$ is bounded with respect to
  $\Vert\cdot\Vert$, then $\phb^{-1}(B)\cap\GZ[1/\beta]$ is locally
  finite. 
\end{lemma}
\begin{proof}
Let $B$ be a bounded subset of $\GK_\beta$, and $x\in\GQ(\beta)$ such that
$\phb(x)\in B$. In particular, for every $i$, $1\leq i\leq\nu$, there
exists a rational integer $m_i$, such that the embedding of $x$ in
$\GK_{\mathfrak{P}_i}$ has valuation at most $m_i$. For $m=\max_{1\leq
  i\leq\nu}m_i$, we get $\beta^m x\in\prod_{v\in\Sb}\mathfrak{O}_v$. On the
other hand, $\beta$ is a $\Sb$-unit, so that
$\beta^m\GZ[1/\beta]=\GZ[1/\beta]\subset\mathfrak{O}_\mathfrak{P}$ for any
$\mathfrak{P}$ coprime with $(\beta)$. Therefore, $\beta^m x\in\mathfrak
O$. Furthermore, the Archimedean absolute values $\vert \beta^m
x\vert_{\beta_i}$ are bounded as well for $i=2,\ldots, r+s$. If we assume
further that $x$ belongs to some bounded subset of $\GQ(\beta)$ (w.r.t. the
usual metric), then all the conjugates of $\beta^m x$ are bounded. Since these
numbers belong to $\mathfrak{O}$, there are only a finite number of them.
\end{proof}
\subsection{Complete tiles and an Iterated Function system}\label{subsec:ct}

\begin{definition}[Complete tiles]\label{def:map}
The complete tiles are the analogues in $\GK_\beta$ of the Euclidean
tiles:
\begin{itemize}
\item Complete central tile
\[
\Tt =  \overline{\phb(\Int(\beta))}
\subset\prod_{v\in\Sb}\OO_v. 
\] 
\item Complete $x$-tiles. For every $x \in {\mathbb Z}[1/\beta]\cap[0,1)$,
\[
\Tt(x) =  \overline{ \phb(\{y \in {\mathbb R}^+;
  \,  \{y\}_\beta=x \}) } \subset \phb(x) +\widetilde{\mathcal
  T}.\; \mbox{In particular, $\widetilde{\mathcal T}=\widetilde{\mathcal
    T}(0)$.}
\] 
\item Complete central subtiles. For every $a \in \{1, \dots, n\}$,
\[
\Ta{a} =\overline{
  \phb\left(\left\{x\in\Int(\beta);\,
  d_\beta(x)\in\La{a}\right\}\right)}.
\]
\end{itemize}
\end{definition}
Using (\ref{T(a)}), we get:
\begin{align}\label{eq:translates}
\Tt(x)&=\phb(x)+\overline{\phb\left(\left\{y \in\Int(\beta);
    d_\beta(y)\cdot d_\beta(x)\in\Lb^{\infty}\right\}\right)}\notag\\
 &=\phb(x)+\overline{\bigcup_{a;\, t_1\cdots t_{a-1}\cdot
  d_\beta(x)\in\Lb^{\infty}}\phb\left(\left\{y\in\Int(\beta);\;
  d_\beta(y)\in\La{a}\right\}\right)}\\ 
&= \phb(x)+\bigcup_{a;\; x<T_\beta^{(a-1)}(1)}\Ta{a}.\notag
\end{align}
Hence, any complete $x$-tile is a finite union of translates of complete
central subtiles.\\
 We now consider  the following  self-similarity property satisfied by the
 complete central subtiles:
\begin{proposition}\label{prop:ifs}
Let $\beta$ be a Pisot number.
The complete  central subtiles satisfy an Iterated Function System equation
(IFS) directed by the admissibility graph (drawn in
Figure~\ref{Fig:admissibility}) in which the direction of edges is
reversed: 
\begin{equation}\label{eq:ifs}
\Ta{a}=\bigcup_{b\xrightarrow{\varepsilon}a}
\left(\phb(\beta)\Ta{b}+\phb(\varepsilon)\right).  
\end{equation}
We use  here Notation~\ref{not:rn}  and we  recall that the digits
$\varepsilon$  belong to  
$\mathcal{D}=\{0,\ldots, \lceil \beta \rceil -1\}$, and that the nodes $a,b$
belong to $\mathcal{A}=\{0,1,\cdots, n\}$.
\end{proposition} 
\begin{proof}
The following decomposition of the languages $\La{a}$ can be read off from the
admissibility graph~\ref{Fig:admissibility} 
\begin{equation}\label{eq:ifslanguages}
\La{a}=\bigcup_{b\xrightarrow{\varepsilon}a}
\La{b}\cdot\{\varepsilon\}.  
\end{equation}
That decomposition yields a similar IFS as in~(\ref{eq:ifs}) where the complete
central subtiles $\Ta{a}$ are replaced by the  images
$\phb(\{x\in\Int(\beta);\, d_\beta(x)\in\La{a}\})$ of the languages
$\La{a}$ into $\GK_\beta$. Lastly, one gets  (\ref{eq:ifs}) by taking the
closure (the unions are finite). It should be noted that this argument does not
depend on the embedding; it is therefore the same as in the unit case, that can
be found  e.g. in \cite{SW02,Siegel03, BertheSiegel05}.
 \end{proof}
\begin{remark} If one details the IFS given by (\ref{eq:ifs}), this gives (with
  $m$ defined in Notation~\ref{not:rn}):
\begin{equation}\label{eq:ifslong}
 \left\{ 
\begin{array}{rcl}
\Ta{1} & = &
\bigcup_{a\in\A}\bigcup_{\varepsilon<t_a}
\left(\phb(\beta)\Ta{a}+\phb(\varepsilon)
 \right) \\ 
\Ta{r+1} & = & \left( \phb(\beta)
  \Ta{m} + \phb(t_{m}) \right) \bigcup\left( \phb(\beta) \Ta{n} +
  \phb(t_{n}) \right) \\ 
\Ta{k+1} & = & \phb(\beta) \Ta{k} +  \phb(t_{k})  ,\;
k\in\{1,\ldots,n-1\}\setminus\{m\}. 
\end{array} \right. 
\end{equation}
\end{remark}

\subsection{Boundary graph}\label{sec:boundarygraph}
The aim of this section is to  introduce the notion of boundary graph which
will be a crucial tool for our estimations of the  function $\gamma$ in
Section~\ref{sec:examples}. This graph is based on the self-similarity
properties  of the boundary of the central tile, in the spirit 
of the those defined  in \cite{Siegel03,Thuswaldner06,SiegelThuswaldner}. 
The idea is the following: in order to understand better the
covering~\eqref{eq:covunit}, we need to exhibit which points belong to the
intersections between the central tile $\Tt$ and the $x$-tiles $\Tt(x)$. To do
this, we first decompose $\Tt$ and $\Tt(x)$ into subtiles: we know 
that $\Tt = \cup_{a \in\A} \Ta{a} $ and Eq.~\eqref{eq:translates} gives $\Tt(x)
= \cup_{b \in {\mathcal A},   T^{b-1}(1)>x} \Ta{b}+\phb(x).$  
Then the intersection between $\Tt$ and $\Tt(x)$ is the union of intersections
between $\Ta{a}$ and $\Ta{b}+\phb(x)$ for  $T^{b-1}(1)>x$. We build a graph
whose nodes stand for each intersection of that type, hence the nodes are
labelled by triplets $[a,x,b]$. To avoid the non-significant intersection
$\Ta{a}\cap\Ta{a}$, we will have to exclude the case $x=0$ and $a=b$.
Then we use the self-similar equation Eq.~\eqref{eq:ifs} to decompose the
intersection $\Ta{a} \cap (\Ta{b}+\phb(x))$ into new intersections of the same
nature (Eq.~\eqref{eq:bigdec}). An edge is labelled with couple of digits, so
that a jump from one node to an another one acts as 
a magnifier of size $\beta^{-1}$, the label of the edge sorting
one digit of the element in the intersection we are describing.

By applying this process, we show below that we obtain a graph that describes
the intersections  $\Ta{a} \cap (\Ta{b}+\phb(x))$
(Theorem~\ref{thm:boundarygraph}). It can be used to check whether  the
covering (\ref{eq:covunit}) is a tiling, as was done in
\cite{Siegel03,SiegelThuswaldner} but this is not the purpose of the present
paper. It the last section, we will use this graph to deduce information on
pure periodic expansions. 

%This graph will be a   very efficient tool   for the study of  the intersection
%between a  complete central subtile  $\Ta{a}$ and  translates of the form 
%$\Ta{b} + \phi_{\beta}(x)$, for $x \in \GZ[1/\beta] \cap [0,1)$; such
%translates  will occur by   extending  the covering (\ref{eq:covunit}) to the
%case of complete tiles. 

\begin{definition}  \label{def:bg}
The nodes of the {\em boundary graph} are the triplets $[a,x,b] \in \A\times
\GZ[1/\beta] \times \A$ such that:
\begin{enumerate}
\item[(N1)] $-T_{\beta}^{(a-1)}(1)<x< T_{\beta}^{(b-1)}(1)$ and $a\not=b$ if $x=0$. 
\item[(N2)] $\phb(x)\in\Ta{a}-\Ta{b}$.
%\item For every conjugate $\beta_i$ of $\beta$, let $\tau_i$ denote the
%  canonical morphism from ${\mathbb Q}(\beta)$ to  ${\mathbb Q}(\beta_i)$. 
%Then $x$ satisfies $$|\tau_i(x)| \leq \frac{[\beta]}{1 - |\beta_i|} $$  
%\item $x \in {\mathcal O}_{\mathbb K}$.
\end{enumerate}
The  labels of the edges  of the boundary graph belong  to $\mathcal{D}^2$.
There exists an edge  $[a,x,b]\xrightarrow{(p_1,q_1)}[a_1,x_1,b_1]$ if and 
only if:
\begin{itemize}
\item[(E1)] $x_1 = \beta^{-1}(x+q_1 - p_1)$,
\item[(E2)] $a_1 \xrightarrow{p_1}a$ and $b_1 \xrightarrow{q_1} b$ are edges of
  the admissibility graph. 
\end{itemize}
\end{definition}

We  first deduce from the definition that the boundary graph is finite  and the
Archimedean norms of its nodes are explicitly bounded: 

\begin{proposition}\label{bounds}
The boundary graph is finite. If $[a,x,b]$ is a node of the  boundary graph,
then we have: 
\begin{itemize}
\item[(N3)] $x \in \OO$;
\item[(N4)] for every conjugate $\beta_i$ of $\beta$,    $|\tau_i(x)|
  \leq \frac{[\beta]}{1 - |\beta_i|}$. 
\end{itemize}
\end{proposition}

\begin{proof}
Let $[a,x,b]$ is a node of the graph. By definition,
$\phb(x)\in\Ta{a}-\Ta{b}$, which  implies 
 $|\tau_i(x)| \leq   \frac{[\beta]}{1 - |\beta_i|}$.

 Let $\mathfrak{P}$ be  a
prime ideal in $\OO$. If  $\mathfrak{P}\mid (\beta)$, then
$x\in\OO_{\mathfrak{P}}$ - since  
$\phb(x)\in\Ta{a}-\Ta{b}$.  Otherwise,   if $\mathfrak{P}$ is coprime
with $\beta$, we use the fact  that
$x\in\GZ[1/\beta]$ to deduce that  $x\in\OO_{\mathfrak{P}}$.  We thus have
$x\in \OO$. It directly follows from Lemma~\ref{lem:finite} that the boundary
graph is finite.
\end{proof}

Proposition~\ref{bounds} will be used in Section~\ref{sec:examples} to
explicitely compute the boundary graph in some specific 
cases: let us stress the fact that  condition (N2) in  Definition~\ref{def:bg}
cannot be directly checked algorithmically, whereas numbers satisfying
condition (N3) and (N4) are explicitely computable. 
Nevertheless, conditions (N3) and (N4) are only necessary conditions for a
triplet to belong to the graph. Theorem~\ref{thm:boundarygraph} below has 
two ambitions: it first details how  the boundary graph indeed describes the
boundary of the graph, as intersections between the central tile and its
neighbours.  Secondly, we will decuce from this theorem an explicit way of
computation for the boundary graph in
Corollary~\ref{cor:boundarygraphalgorithm}.
 
The following lemma shows that  Condition (N1)  in  Definition~\ref{def:bg}
 automatically  holds    for a node $[a_1,x_1,b_1] \in \A\times
\GZ[1/\beta] \times \A$ as soon as one has  the  edge conditions  between
$[a,x,b] \in \A\times \GZ[1/\beta] \times \A$ and $[a_1,x_1,b_1]$.

\begin{lemma}\label{lem:graphbound}
Let $x \in \left(-T_{\beta}^{(a-1)}(1), T_{\beta}^{(b-1)}(1)\right)
\cap\GZ[1/\beta]$. Let $a_1 \xrightarrow{p_1} a$ and $b_1 \xrightarrow{q_1} b$  
be two edges in the  admissibility graph.  Let  $x_1 = \beta^{-1}(x+q_1 -
p_1)$. One has $x_1 \in \left(-T_{\beta}^{(a_1-1)}(1),
  T_{\beta}^{(b_1-1)}(1)\right)$.  
\end{lemma}
\begin{proof}
Assume that $x$ is non-negative (otherwise, the same argument
applies to $-x$).  We thus have  $- \frac{p_1}{\beta}  \leq x_1 \leq
\frac{x+q_1}{\beta}$. 
Since $a_1\xrightarrow{p_1}a$, we have that $p_1\leq
t_{a_1}$, hence $p_1\, 0^{\infty}<_{\rm lex}
S^{a_1-1}(d_\beta^*(1))$ (the strict inequality comes from the fact
that  $d_\beta^*(1))$ does not  ultimately end  in  $0^{\infty}$).  
 Therefore, $x_1\geq
-\frac{p_1}{\beta}>-T_{\beta}^{(a_1-1)}(1)$ by~(\ref{T(a)}).  

On the other hand, since $x <  T_{\beta}^{(b-1)}(1)$, then
the sequence  $t_1\cdots t_{b-1}\, d_\beta(x)$ 
is admissible, again by~(\ref{T(a)}).   We  thus deduce from
 $b_1\xrightarrow{q_1}b$  that  $t_1\cdots t_{b_1-1}\,(q_1\,
d_\beta(x))$  is admissible. We thus get
$x_1\leq\frac{x+q_1}{\beta}<T_{\beta}^{(b_1-1)}(1)$. 
\end{proof}
However, if $\beta$ is not a unit, it does not follow from
Lemma~\ref{lem:graphbound} that if $[a,x,b]$ 
is a node of the boundary graph, $a_1 \xrightarrow{p_1}a$ and $b_1
\xrightarrow{q_1} b$ are edges of the admissibility graph, and $x_1 =
\beta^{-1}(x+q_1 - p_1)$, then $[a_1,x_1,b_1]$ is a node (we have also  to
check  Condition  (N2) or (N3)): for instance, consider the two edges of the
admissibility graph $1\xrightarrow{0}1$ and $1\xrightarrow{t_1}2$. Starting
from the note $[1,0,2]$, the edges above would yield
$x_1=-\frac{t_1}{\beta}\not\in\OO$. Hence
$[1,x_1,1]$ is not a node of the boundary graph by Proposition~\ref{bounds}.\\

We now prove that the boundary graph is indeed a good description of the
boundary of the central tile,  by relating it with intersections  between
translates of the  complete  central subtiles. 
\begin{theorem}\label{thm:boundarygraph}
%Let $x \in (-T_{\beta}^{(a-1)}(1), T_{\beta}^{(b-1)}(1)) \cap\GZ[1/\beta]$.
 Let $z \in\GK_{\beta}$.  The point $z$ belongs 
to the intersection $\Ta{a}\cap(\Ta{b}+\phb(x))$, for $x\in \GZ[1/\beta]$,
with $a\not=b$ if $x=0$, if and only if  $[a,x,b]$ is a node of the graph and
there exists an infinite path in the boundary graph, starting from the node  
$[a,x,b]$ and labeled by $(p_i,q_i)_{i\geq 0}$ such that $$z = \sum_{i=0}^\infty
\phb(p_i \beta^i).$$
\end{theorem}

\begin{proof}
Let $x \in (-T_{\beta}^{(a-1)}(1), T_{\beta}^{(b-1)}(1)) \cap\GZ[1/\beta]$.
The  complete central subtiles satisfy a  graph-directed self-affine equation
detailed in Proposition~\ref{prop:ifs} that yields the decomposition
\begin{equation}\label{eq:bigdec}
\Ta{a}\cap\left(\Ta{b}+\phb(x)\right)=
\bigcup_{\genfrac{}{}{0pt}{1}{a_1\xrightarrow{p_1}a}{b_1\xrightarrow{q_1}b}} 
\left[
\left(\phb(\beta)\Ta{a_1}+\phb(p_1)\right)\cap
\left(\phb(\beta)\Ta{b_1}+\phb(q_1)+\phb(x)\right)
\right].
\end{equation}
Let $z \in \Ta{a}\cap(\Ta{b}+\phb(x))$. Then there exist two edges
$a_1\xrightarrow{p_1}a$ and  $b_1\xrightarrow{q_1} b$ such that the
corresponding intersection in the right-hand side of~(\ref{eq:bigdec}) contains
$z$. Setting $x_1=\beta^{-1}(x+q_1-p_1)$ and
$z_1=\phb(\beta)^{-1}(z-\phb(p_1))$, we get
$z_1\in\Ta{a_1}\cap\left(\Ta{b_1}+\phb(x_1)\right)$. By construction,
$x_1\in\GZ[1/\beta]$ and belongs to the interval  $(-T_{\beta}^{(a_1-1)}(1),
T_{\beta}^{(b_1-1)}(1))$ by Lemma~\ref{lem:graphbound}. Then, by definition,
$[a_1,x_1,b_1]$ is a node of the boundary graph, and we may iterate the above
procedure. After $n$ steps, we have
\[
\frac{z-\phb\left(\sum_{i=1}^np_i\beta^{i-1}\right)}
{\phb(\beta^n)}\in \Ta{a_n}\cap\left(\Ta{b_n}+\phb(x_n)\right).
\]
It follows that
$\Vert z-\phb\left(\sum_{i=1}^np_i\beta^{i-1}\right)\Vert
\ll\Vert\phb(\beta)\Vert^n$ for $n$ tending to infinity; therefore
$z=\sum_{i=1}^\infty\phb(p_i\beta^{i-1})$.\\

Conversely, let $z$ such that  $z= \sum_{i\geq 1} \phb(\beta^{i-1}
p_i)$ with $(p_i,q_i)_{i \geq 1}$ the labeling of a path on the
boundary graph starting from $[a,x,b]$. By the definition of the edges of the
graph,  one checks  that  $t_1\cdots t_{a-1}$ 
is a suffix of $t_1\cdots t_{a_1-1}p_1$, which is itself suffix of $t_1\cdots
t_{a_2-1}p_2p_1$, and so on. Hence $z\in\Ta{a}$. Let $y=\sum_{i\geq 1}
\phb(\beta^{i-1}q_i)$. By construction, we also have $y\in\Ta{b}$. Furthermore,
the recursive definition of the $x_i$'s gives
$$x+\sum_{i=1}^nq_i\beta^{i-1}=\sum_{i=1}^np_i\beta^{i-1}+\beta^nx_n.$$
The sequence $(x_n)_n$ takes only finitely many values by Proposition
\ref{bounds}, hence $\phb(\beta^nx_n)$ 
tends to 0, which yields $\phi_\beta(x)+y=z$. Therefore
$z\in\Ta{a}\cap\left(\Ta{b}+\phi_\beta(x)\right)$. 
\end{proof}

\begin{corollary} \label{cor:inter}
Let $x \in \GZ[1/\beta]$ and $a\neq b$ if $x=0$.   The intersection
$\Ta{a}\cap(\Ta{b}+\phb(x))$ is 
non-empty  if and only  if $[a,x,b]$  is a node 
of the boudary graph and there exists at least an  infinite path in  the
boudary graph starting from  $[a,x,b]$.
\end{corollary}

We deduce a procedure for the computation of the boundary graph.

\begin{corollary} \label{cor:boundarygraphalgorithm}
The boundary graph can be obtained as follows:
\begin{itemize}
\item Compute the set of triplets $[a,x,b]$ that satisfy conditions (N1), (N3)
  and (N4); 
\item Put edges between two triplets if conditions (E1) and (E2) are satisfied;
\item Recursively remove nodes that have no outging edges.
\end{itemize}
\end{corollary}

%and both properties proved in  Proposition~\ref{bounds} and whose edges
%correspond to Definition~\ref{def:bg}. 
\begin{proof}
The particularity of this graph is that any node belongs to an infinite path.
Proposition~\ref{bounds} and  Theorem~\ref{thm:boundarygraph} show that this
graph is bigger than (or equal to) the boundary graph. Nevertheless, 
the converse part of the proof of Theorem~\ref{thm:boundarygraph} ensures that
if an infinite path of the latter graph starts from $[a,x,b]$, then this path
produces an element $z$ in $\Ta{a}\cap(\Ta{b}+\phb(x))$. Therefore,
$\phb(x)\in\Ta{a}-\Ta{b}$, and $[a,x,b]$ is indeed a node of the boundary
graph. Finally, even if the procedure described in the statement of the
corollary  mentions infinite paths, it needs only finitely many operations,
since the number of nodes is finite: it has been proved in
Proposition~\ref{bounds} for the boundary graph; it is an immediate consequence
of Lemma~\ref{lem:finite} for triplets satisfying (N1), (N3) and (N4).
\end{proof}

\subsection{Covering of the complete representation space}

In order to generalise the tiling property stated in Theorem~\ref{thm:Pavageunit} to the non-unit case, we
need to understand better what is the complete representation of $\GZ[1/\beta] 
\cap \RR^+$. We first prove the following lemma, that makes
Lemma~\ref{lem:density} more precise.

\begin{lemma}\label{lem:covering}
We have that $\phb\left(\mathfrak{O} \cap\GR^{+}\right)$ is dense in
$\prod_{v\in\Sb}\mathfrak{O}_v$ and that $\phb\left(\GZ[1/\beta]
  \cap\GR^{+}\right)$ is dense in $\GK_{\beta}$. Those density results remain
true if one replaces $\GR_+$ by any neighbourhood of $+\infty$.  
\end{lemma}

\begin{proof}
We already know by Lemma~\ref{lem:density} that $\phb(\mathfrak{O})$ is
dense in $\prod_{v\in\Sb} \mathfrak{O}_v$. Let $U\geq 0$. For any
$x\in\mathfrak{O}$, we have $x+\beta^n> U$ if $n$ is sufficiently large. Since
$\beta^n$ tends to 0 in $\GK_\beta$, $\phb(x+\beta^n)$ tends to
$\phi(x)$; hence $\phb\left(\mathfrak{O} \cap[U,+\infty)\right)$ is dense in
$\prod_{v\in\Sb}\mathfrak{O}_v$.

Let  $Z=(z,y_1, \dots, y_{\nu}) \in {\mathbb K}_{\beta}$. Since ${\mathbb
  K}_{\beta}$ is built from the prime divisors of $\beta$, there exists a
natural integer $n$ such that
$\beta^n y_i \in \mathfrak{O}_{\mathfrak{P}_i}$ for every $i=1,...,\nu$.
Moreover, there exists an integer $A$ such that  
$A\mathfrak{O}\subset {\mathbb Z}[\beta]$ (for instance, the discriminant of
$(1,\beta,\ldots,\beta^{d-1})$). Split $A$ 
into $A=A_1 A_2$, so that $A_1$ is coprime with $\beta$ and 
the prime divisors of $A_2$ are also divisors of $N(\beta)$.  Then $A_1$ is a
unit in each $\mathfrak{O}_{\mathfrak{P}_i}$
so that $y_i/A_1 \in\mathfrak{O}_{\mathfrak{P}_i} $ for $1\leq i\leq\nu$.  By
the definition of $A_2$, there exists $m$ such that $\beta^m/A_2
\in\mathfrak{O}$. Therefore, $\beta^{\max(n,m)}Z/A_2\in\prod_{v\in\Sb}
\mathfrak{O}_v$. Applying the first part of the lemma, there exists a sequence
$(x_l)_l$ in $\GZ[\beta^{-1}]\cap[U\beta^{\max(n,m)},+\infty)$ such that
$(\phb(x_l))_l$ tends to $\beta^{\max(n,m)}Z/A_2$. Then,
$(\phb(\beta^{-\max(n,m)}A_2x_m))_l$ tends to $Z$. Since
$\beta^{-\max(n,m)}A_2x_l\in\GZ[1/\beta]\cap[U,+\infty)$, the proof is complete. 
\end{proof}

\begin{proposition}\label{prop:covering}
The complete central tile $\Tt$ is compact. The $x$-tiles $\Tt(x)$ provide a
covering of the $\beta$-representation space: 
\begin{equation}\label{eq:tiling} 
\bigcup_{x \in {\mathbb Z}[1/\beta]\cap[0,1)}\Tt(x) =\GK_{\beta}.
\end{equation}
Moreover, this covering is uniformly locally finite: for any $R>0$, there
exists $\kappa(R)\in\GR_+$ such that, for all $z\in\GK_\beta$, one has
\[
\#\left\{x \in\GZ[1/\beta]\cap[0,1);\; \Tt(x)\cap
B(z,R)\neq\emptyset\right\}\leq\kappa(R). 
\]
\end{proposition}

\begin{proof}
The projection of $\Tt$ on $\GK_f$ is compact since the local rings $\OO_v$
are. Its projection on $\GK_\infty$ is bounded because $\beta$ is a Pisot
number. Since $\Tt$ is obviously closed, it is therefore compact. Explicitly, 
we have by construction that $\Vert\phb(\beta)\Vert<1$. Since $\Vert
n\Vert=n$ for each $n\in\GZ$, it follows that $\Tt\subset B(0,M_1)$ with $M_1
=(\lfloor\beta\rfloor)/(1-\Vert\phb(\beta)\Vert)$. 

Since $\beta$ is an integer, we have
$\Int(\beta)\subset\GZ[1/\beta]$. Therefore, for $y\in\GR_+$, $y$ belongs to
$\GZ[1/\beta]$ if and only if $\{y\}_\beta$ belongs to $\GZ[1/\beta]$. In other
words, 
\[
\bigcup_{x\in\GZ[1/\beta]\cap[0,1)}\big\{y\in\GR_+;\, \{y\}_\beta
=x\big\}=\GZ[1/\beta]\cap\GR_+,
\]
and, by Lemma~\ref{lem:covering}, we have that
\begin{equation}\label{eq:finiteunion}
{\mathbb K}_{\beta} = \overline{\phb({\mathbb Z}{[1/\beta]\cap
    {\mathbb R}^+) }} =  
\overline{\bigcup_{x \in {\mathbb Z}[1/\beta]\cap[0,1)}
  \phb\left(\big\{y\in\GR_+;\, \{y\}_\beta =x\big\} \right) }. 
\end{equation}
Let us fix $z \in \GK_{\beta}$ and $R>0$. We consider the ball
$B(z,R)$ in $\GK_{\beta}$.
Assume  that $x\in\GZ[1/\beta]\cap[0,1)]$ is such that $\Tt(x)\cap
B(z,R)\neq\emptyset$. By
$\Tt(x)\subset\phb(x)+\Tt\subset\phb(x)+B(0,M_1)$. Hence  
$\phb(x)\in B(z,R+M_1)$. Then, Lemma~\ref{lem:finite} ensures that there
exists only finitely many such $x$. 

It certainly remains to prove that the number of those $x$ is bounded
independently of $z$, but it already shows that the union in the right-hand
side of~(\ref{eq:finiteunion}) is finite, which allows to permute the union and
the closure operations and proves~(\ref{eq:tiling}). 

We then
use~(\ref{eq:tiling}) to prove the existence of some
$x_0\in\GZ[1/\beta]\cap[0,1]$ such that
$\Vert\phi_\beta(x_0)-z\Vert<1$. Therefore, any $x\in\GZ[1/\beta]\cap[0,1)$
satisfying $\Tt(x)\cap B(z,R)\neq\emptyset$ can be written as $x=x_0+x_1$,
where $x_1\in\GZ[1/\beta]\cap[-2,1)$ and $\phb(x_1)\in
B(0,R+M_1+1)$. Lemma~\ref{lem:finite} gives an upper bound $\kappa(R)$ for the
number of such $x_1$, and the lemma is proved.
\end{proof}

\begin{corollary}\label{cor:covering}
The complete central tile $\Tt$ has non-empty interior in the representation
space ${\mathbb K}_{\beta}$, hence non-zero Haar measure. 
\end{corollary} 
\begin{proof}The property concerning the complete central tile has already been
  proved in \cite{BertheSiegel07}, Theorem~2-(2), by geometrical
  considerations. However, most of this proposition is now an immediate
  consequence of~(\ref{eq:tiling}): since $\GK_\beta$ is locally compact, it is
  a Baire space. Therefore, some $\Tt(x)$ must have non-empty interior, hence
  the central tile itself, by $\Tt(x)\subset\phb(x)+\Tt$. Thus it has
  positive measure. By the way, (\ref{eq:tiling}) gives also a direct proof of
  that fact without any topological consideration, by using the
  $\sigma$-additivity of the measure and $\mu_\beta(\Tt(x))\leq \mu_\beta(\Tt)$.
 \end{proof}

\subsection{Inner points}

We use the covering property to express the complete central tile as the
closure of its exclusive inner points (see Definition
\ref{def:exclu}). Since we will use it extensively, we
introduce the notation $c_\beta=\Vert\phb(\beta)\Vert.$ We have seen that
$0<c_\beta<1$. 

\begin{proposition}\label{prop:innerpoint}
Let $\beta$ be a Pisot number.
If $\beta$ satisfies the property (F), then $0$ is an exclusive inner point of
the complete central tile $\Tt$. Indeed, it is an inner point of the complete central
subtile $\Ta{1}$.
\end{proposition}
\begin{proof}
By Lemma~\ref{lem:finite}, there exist  finitely many
$x\in\GZ[1/\beta]\cap[0,1)$ such that $\Vert\phb(x)\Vert\leq
2M_1$, where the constant $M_1$ is taken from the proof of
Proposition~\ref{prop:covering}. According to property (F), all those $x$ have
finite $\beta$-expansion. Let $p$ be the maximal length of those expansions. 

Let $m$ be a  non-negative integer and $x\in(\GZ[1/\beta]\cap\GR_+)\setminus
\beta^m\Int(\beta)$. Set 
$x_1=\lfloor\beta^{-p-m}x\rfloor_\beta$ and $x_2=\{\beta^{-p-m}x\}_\beta$. By
construction, we have $\Vert\phb(x_1)\Vert\leq M_1$ and
$\Vert\phb(x_2)\Vert>2M_1$, the latter because $d_\beta(x_2)$ has length
greater than $p$.  Set $M_2= M_1 c_{\beta} ^p $.Therefore, we have that
\[
\Vert\phb(x)\Vert=
c_\beta^{p+m}\Vert\phb(x_1)+\phb(x_2)\Vert>
M_1c_\beta^{p+m}=c_\beta^mM_2.
\]
Hence, we have $\phb^{-1}\big(B(0,c_\beta^mM_2)\big)\cap\GZ[1/\beta]\cap\GR_+
\subset\beta^m\Int(\beta)$. Taking $m=0$, this shows that the origin is
exclusive. Moreover, since $B(0,c_\beta^m M_2)$ is open, and since
$\phi_\beta(\beta^m\Int(\beta))\subset\Ta{1}$ for $m$ sufficiently large,
Lemma~\ref{lem:covering} ensures that $B(0,c_\beta^m M_2)\subset\Ta{1}$. 
\end{proof}

\begin{theorem}\label{thm:disjointness}
Let $\beta$ be a Pisot number. Assume that $\beta$ satisfies the finiteness
property (F). Then each tile $\Tt(x)$, $x \in\GZ[1/\beta]\cap[0,1)$ 
is the closure of its interior, and each inner point of  $\Tt(x)$ is
exclusive. Hence, for every $x\not=x^{\prime} \in\GZ[1/\beta]\cap[0,1)$,
$\Tt(x')$ does not intersect the interior of  $\Tt(x)$. The
tiles $\Tt(x)$ are measurably disjoint in $\GK_{\beta}$. Moreover, their
boundary has zero  measure. 

The same properties hold for  the translates of  complete central subtiles
$\Ta{a}+ \phi_{\beta} (x)$, for $a\in\A$ and $x \in \GZ[1/\beta] \cap [0,1)$.
\end{theorem}
\begin{proof}
The proof of the unit case can be found in \cite{Akiyama02}(Theorem 2,
Corollary 1) and could have been adapted. We follow here a slightly different
approach. For $x\in\GZ[1/\beta]\cap\GR_+$, let 
$Y(x)=\{y\in\GR_+;\, \{y\}_\beta=x\}\subset\GZ[1/\beta]\cap\GR_+$. By
definition, $\Tt(x)=\overline{\phb(Y(x))}$. According to the proof of
Proposition~\ref{prop:innerpoint}, we have
\begin{equation}\label{eq:lesboules}
{\phb}^{-1}(B(0,c_\beta^m
M_2))\cap\GZ[1/\beta]\cap\GR_+ \subset\beta^m\Int(\beta)\quad\&\quad
B(0,c_\beta^m 
M_2)\subset\overline{\phb(\beta^m\Int(\beta))}. 
\end{equation}
Recall that $n$ is the length of $d_\beta(1)$. Therefore, if $w_1$ and $w_2$
are admissible, so is $w_1\cdot 0^n\cdot w_2$. Now, for any given $y\in Y(x)$,
we have $y+\beta^{m}\Int(\beta)\subset Y(x)$ for any 
$m\geq m(y)=n+\lceil(\log y)/(\log\beta)\rceil$. Therefore,
\begin{equation} \label{eq:39bis}
\Tt(x)=\overline{\bigcup_{y\in
    Y(x)}\phb\big(y+\beta^{m(y)}\Int(\beta)\big)} =\overline{\bigcup_{y\in
    Y(x)}B\big(\phb(y),c_\beta^{m(y)}M_2\big)},
\end{equation}
and $\Tt(x)$ is the closure of an open set, hence of its
interior.\\

Therefore, in order to prove the exclusivity, we only have to show that two
different $x$-tiles have disjoint interiors. Let
$x,x'$ in $\GZ[1/\beta]\cap[0,1)$. According to (\ref{eq:39bis}), 
any non-empty open subset of $\Tt(x)\cap\Tt(x')$ contains some 
ball $B=B(\phb(y), c_\beta^m M_2)$, with $y\in Y(x)$ and $m=m(y)$ chosen as
above. Since $\phb$ is a ring homomorphism, the first part
of~(\ref{eq:lesboules}) implies that
${\phb}^{-1}(B)\cap\GZ[1/\beta]\cap[y,+\infty) \subset
y+\beta^m\Int(\beta)$. But there also exists $y'\in Y(x')$ and a natural 
integer $m'$ such that $\phb(y'+\beta^{m'}\Int(\beta))\subset B$. Since
$y'+\beta^{m'}\Int(\beta)$ contains arbitrary large real 
numbers, this shows that $Y(x)\cap Y(x')\neq\emptyset$. Hence $x=x'$ and the
exclusivity follows.\\

The proof for the subtiles $\widetilde{\mathcal T}^{(a)}$ works exactly in the
same way, because of the key property 
\[
d_\beta(x)\in\La{a}\Longrightarrow\forall y\in\Int(\beta):\,
d_\beta(x+\beta^my)\in\La{a}
\]
for $m$ sufficiently large (depending on $x$).\\

It is possible to prove directly that the subtiles $\Ta{a}$ are measurably
disjoint (for an efficient proof based on the IFS~(\ref{eq:ifs})  and
Perron-Frobenius Theorem, see \cite{SW02,BertheSiegel05}[Theorem~2]). However,
it follows directly from the fact that the boundary of the subtiles have
zero-measure, since two different subtiles have disjoint
interiors.

To prove the latter, we follow~\cite{Praggastis} [Proposition~1.1]. Since $\A$
is finite, there exist $\delta$ and $a\in\A$ such that 
$\mu_\beta(\partial\Ta{a})=\delta\mu_\beta(\Ta{a})$ and
$\mu_\beta(\partial\Ta{b})\leq\delta\mu_\beta(\Ta{b})$ for all $b\in\A$.
Let $k\geq n$ be a rational integer. Then, by~(\ref{eq:translates}), we have
\begin{align}\label{eq:decomposition}
\phb(\beta)^{-k}\Ta{a}&=\overline{\big\{\phb(\beta^{-k}x);\, x\in\Int(\beta),
d_\beta(x)\in\La{1}\big\}}\notag\\
&=\bigcup_{x\in\Lambda_k}\Tt(x)
=\bigcup_{x\in\Lambda_k} \bigcup_{b;\,
x<T^{(b-1)}(1)}\left(\phb(x)+\Ta{b}\right),
\end{align}
where $\displaystyle\Lambda_k=\left\{\sum_{i=1}^k\omega_i\beta^{-i};\,
\omega_1\cdots\omega_k\in\La{a} \right\}$. The $x$-tiles (resp. the subtiles)
having disjoint interiors, the family of tiles $\phb(x)+\Ta{b}$ occurring
in~(\ref{eq:decomposition}) has the same property. Then, for a subfamily
$(\mathcal{T}_i)_i$ of those tiles, we have $\mathcal{T}_i\cap
\mathcal{T}_j=\partial \mathcal{T}_i\cap\partial\mathcal{T}_j$, 
and a simple argument gives
$\mu_\beta(\partial(\cup\mathcal{T}_i))\leq\delta\mu_\beta(\cup
\mathcal{T}_i)$. Let us split the 
union~(\ref{eq:decomposition}) as $\phb(\beta)^{-k}\Ta{a}=\mathcal{U}_1\cup
\mathcal{U}_2$, where $\mathcal{U}_1$ is the union of those tiles intersecting
the boundary of $\phb(\beta)^{-k}\Ta{a}$ and $\mathcal{U}_2$ the union of those
tiles included in its interior. If $k$ is large, $\phb(\beta)^{-k}\Ta{a}$
contains open balls of sufficiently large size to contain some of the tiles,
whose diameter are at most $\max_{b\in\A}{\rm diam}(\Ta{b})$. Hence
$\mathcal{U}_2$ is not empty, and has actually positive measure. Finally, since
the multiplication by $\phb(\beta)$ preserves the boundary, we have
\[
\delta
\mu_\beta\left(\phb(\beta)^{-k}\Ta{a}\right)=
\mu_\beta\left(\partial(\phb(\beta)^{-k}\Ta{a})\right)
\leq\mu_\beta(\mathcal{U}_2)<\delta
\mu_\beta\left(\phb(\beta)^{-k}\Ta{a}\right),
\]
if $\delta\neq 0$, which would yield a contradiction. The metric disjointness
follows for the tiles $\Ta{a}$, hence for the $\Tt(x)$ too by
(\ref{eq:translates}). 
\end{proof}

We can project this relation on the Euclidean space.

\begin{corollary}
Let $\beta$ be  a Pisot number. If $\beta$ satisfies the  finiteness property
(F), then $0$ is an inner point of the central tile ${\mathcal T}$ and each
tile ${\mathcal T}(x)$ is the closure of its interior.  
\end{corollary}

\begin{proof}
If 0 in an inner point of $\widetilde{\mathcal T}$ in the field ${\mathbb
  K}_{\beta}$, then 0 is also an inner point in its projection on ${\mathbb
  K}_{\infty}$.
\end{proof}

This corollary is the most extended generalisation of  Theorem \ref{thm:Pavageunit}
to the non-unit case: if we only consider Archimedean embeddings to build the central tile, 
the finiteness property still implies that  $0$ is  an inner point of the central tile. Nevertheless, 
inner points are no more exclusive, hence the tiling property is not satisfied. 

\medskip

\paragraph{\bf Choosing the suitable non-Archimedean embedding} We already
explained that the Archim\-ed\-ean 
embedding was not suitable for building a  measure-preserving algebraic
extension. We shall comment now 
why the choice of  the beta-adic representation space $\GK_\beta$ 
 is suitable from the tiling viewpoint. 
It is a general fact read only from the admissibility
  graph that the (complete) subtiles $\widetilde{\mathcal T}^{(a)}$ satisfy
  an Iterated Function System (IFS). Thanks to  the introduction of
  the beta-representation space, the action of the multiplication by
  $\phb(\beta)$ in $\GK_{\beta}$ 
acts on the measure as a multiplication  by a ratio $1/\beta$ 
according to~(\ref{equ:homothetie}). That property allows to deduce from the
IFS that the (complete) subtiles are measurably disjoint in
${\mathbb K}_{\beta}$  - whereas their projection  
${\mathcal T}^{(a)}$ on ${\mathbb K}_{\infty}$ are not
(Theorem~\ref{thm:disjointness} below). More geometrically, the space
$\GK_\beta$ is chosen so that:
\begin{itemize}
\item the tiles are big enough to cover it (covering property, first part
of Proposition~\ref{prop:covering}), and  
\item they are small enough, so 
that they do not overlap much - neither combinatorically (locally finitely many
overlaps - second part of Proposition~\ref{prop:covering}) nor
topologically (disjoint interiors), nor metrically (measurable disjointness) -
as shows Theorem~\ref{thm:disjointness} (tiling property).
\end{itemize}
If $\beta$ is not a unit, the space
$\GK_\infty$ is too small to ensure the tiling property. On the opposite, the
restricted topological product of the $\GK_v$ with respect for the $\OO_v$ for
all places $v$ but the Archimedean one given by the identity
embedding (in other words, the projection of the ad\`ele group
$\mathbb{A}_{\GQ(\beta)}$ obtained by
canceling the coordinate corresponding to that Archimedean valuation) would
have satisfied the tiling property and given an interesting algebraical
framework, but would have been too big for the covering property - since the
principal ad\`eles build a discrete subset in the ad\`ele group.

\section{Purely periodic expansions}\label{sec:periodicity}

The elements with a purely periodic expansion, denoted by $\Pi_\beta$ (see Notation \ref{def:periodiques}),  belong to $\GQ(\beta)$ and as explained in the
introduction, 
there are numbers $\beta$ for which $\Pi_\beta^{(r)}=[0,1)\cap \GQ$. However,
Lemma~\ref{lem:pp} below shows that if $\beta$ is a Pisot number, but not a
unit, there  exist arbitrary small rational numbers that do not belong to
$\Pi_\beta^{(r)}$. This justifies the restriction in the definition of 
$\gamma(\beta)$, that only takes into account rational
numbers whose denominator is coprime with the norm of $\beta$. 
\begin{lemma} \label{lem:pp}
Let $\beta$ be a non-unit Pisot number.  Let $x=\frac ab\in\GQ\cap[0,1)$ with
$\gcd(b,N(\beta))>1$. Then $d_\beta(x)$ is not purely periodic.  
\end{lemma}
\begin{proof}
Suppose that the $\beta$-expansion of $x\in\GQ\cap[0,1)$ is purely periodic
with period $l$. Then we can write: 
\[
x=\frac ab= \sum_{k\geq 0}\beta^{-k\ell}(a_1 \beta^{-1} + \cdots +a_\ell
\beta^{-\ell})  = \frac{a_1 \beta^{\ell-1} + \cdots+ a_\ell}{\beta^\ell - 1}.
\]  
Hence $x=\frac{A}{\beta^\ell -1}$ with $A\in\OO$. Since the principal ideals
$(\beta)$ and $(\beta^\ell-1)$ are coprime, we get $\phb(x)\in\prod_v\OO_v$. On
the other hand, if $p\mid \gcd(b, N(\beta)$, then $\phi_\beta(a/b)$ contains a
component in $\GQ_p\setminus\GZ_p$. Hence $a/b\neq A/(\beta^\ell-1)$.
\end{proof}

\subsection{Pure periodicity and   complete tiles}
Using and adapting ideas from~\cite{Praggastis, ItoRao04, Sano02},   one
obtains the following  characterisation of real numbers having a purely
periodic $\beta$-expansion; this result can be considered as a first step
towards the realisation of an algebraic natural extension of the
$\beta$-transformation. Notice that Theorem~\ref{thm:purpernonunit} is
naturally stated in~\cite{BertheSiegel07} with compact intervals, which obliges
to take in account the periodic points and to distinguish whenever $d_\beta(1)$
is finite or infinite. Our point of view simplifies the proof; for that reason,
we give it. 

\begin{theorem}[\cite{BertheSiegel07}, Theorem~3]\label{thm:purpernonunit} 
 Let $x \in [0,1)$. Then, $x$ belongs to $\Pi_\beta$ if  and only if  
\begin{equation}\label{eq:conditionperiodicity}
(-\phb(x),x)
\in \bigcup_{a \in \A}
(\Ta{a})\times\big[0,T_{\beta}^{a-1}(1)\big).
\end{equation}
\end{theorem}
\begin{proof}Let $x \in [0,1)\Pi_\beta$ with purely periodic beta-expansion
  $d_\beta(x)=(a_1\cdots a_\ell)^{\omega}$. Obviously, $x\in\GQ(\beta)$. A
  geometric summation gives 
\[
x=\frac 1{1-\beta^{-\ell}}\sum_{k=1}^\ell
a_k\beta^{-k}=-\frac 1{1-\beta^\ell}\sum_{j=0}^{\ell-1}a_{\ell-j}\beta^j.
\]
Applying $-\phb$ to the latter and going the geometric summation backwards
yields 
\begin{equation}\label{eq:reversion}
-\phb(x)=\sum_{j=0}^\infty\tilde
a_j\phb(\beta)^j=\lim_{n\to\infty}\phb\left(\sum_{j=0}^n\tilde
a_j\beta^j\right),\;\text{with}\; \tilde a_j=a_{\ell-j\pmod\ell}. 
\end{equation}
It is obvious that the sum $\sum_{j=0}^n\tilde a_j\beta^j$ is a beta-expansion,
since we have by construction $\tilde a_{\ell-1}\cdots\tilde a_0=a_1\cdots
a_\ell$. Therefore, $-\phb(x)\in\Tt$. Moreover, the admissibility of the
concatenation $\tilde a_{n}\cdots\tilde a_0\, d_\beta(x)=\tilde
a_{n}\cdots\tilde a_0\tilde a_{-1}\tilde a_{-2}\cdots$ is the exact translation
of the condition $(-\phb(x),x)\in\Ta{a}\times[0,T^{a-1}(1))$. Hence the
condition is necessary.

Let us prove that the condition is sufficient, and let
$z\in\GQ(\beta)\cap[0,1)$ such 
that $(-\phb(z),z)\in\Ta{a}\times [0,T^{a-1}(1))$ for some
$a\in\A$. By compactness, there exists a sequence of digits $(w_n)_n$
such that 
$\phb(z)=\lim_{n\to\infty}\phb(\sum_{j=0}^nw_j\beta^j)$, the latter
sums being beta-expansions for all $n$. Moreover, the bi-infinite word $\cdots
w_nw_{n-1}\cdots w_0\cdot d_\beta(z)$ is admissible.
Define a sequence $(z_k)_k$ by $d_\beta(z_k)=w_{k-1}w_{k-2}\cdots w_0\cdot
d_\beta(z)$. Write $z_0=z=a/b$, with $b\in\GN^*$ and $a\in\GZ[1/\beta]$. Then, 
\begin{equation}\label{eq:zk}
z_k=\beta^{-k}\left(z+\sum_{j=0}^{k-1}w_j\beta^j\right)\in b^{-1}\GZ[1/\beta]. 
\end{equation}
Applying $-\phb$ to~(\ref{eq:zk}) gives
$\displaystyle
{-\phb(z_k)=
\lim_{n\to\infty}\phb\left(\sum_{j=0}^nw_{k+j}\beta^j\right)}$. In
particular, $-\phb(z_k)\in\Tt$ for any $k$, which ensures that the
sequence $(\phb(z_k))_k$ is bounded too. So is the sequence
$(\phb(bz_k))_k$, which is hence finite by Lemma~\ref{lem:finite}. Thus
$z_j=z_{j+s}$ for some $j$ and $s\neq 0$. This shows that
$d_\beta(z)=(w_{s-1}w_{s-2}\cdots w_0)^\infty$ and concludes the proof.
\end{proof}

As shows Theorem~\cite{BertheSiegel07}, the points of the orbit of 1 under the
action of $T_\beta$ play a special role. They have to be treated separately.
\begin{lemma}\label{lem:Torbit}
We have either $T_\beta^k(1)=0$, or $T_\beta^k(1)\in\GQ(\beta)\setminus\GQ$,
but if it is $0$. Moreover,  $T_\beta^k(1)\in\Pi_\beta$ if and only if $\beta$
is a non-simple Parry number (that is $m\neq 0$) and $k\geq m$.
\end{lemma}
\begin{proof} The transformation $T_\beta$ preserves $\OO$. Hence
  $T_\beta^k(1)\in\OO$ for all $k$. Since $\GQ$ is integrally closed, if
  $T_\beta^k(1)\in\GQ$, then $T_\beta^k(1)\in\GZ$. HEnce he only possibility is
  $T_\beta^k(1)=0$. This happens exactly if $\beta$ is a simple Parry number
  (that is if $m=0$) and $k\geq n$.  
We have mentionned in Section~\ref{sec:ag} that
$d_\beta(T_\beta^k(1))=S^k(d_\beta(1))$. Therefore,
$T_\beta^k(1)\in\Pi_\beta$ if and only if $\beta$ is a non-simple Parry
number and $k\geq m$. According to
$d_\beta^*(1)=d_\beta(1)=(t_1\cdots t_m)(t_{m+1}\cdots t_n)^\infty$, the
orbit possesses $n$ elements, $m$ of them having purely periodic beta-expansion.
\end{proof}

\medskip

\paragraph{\bf  Application to the function $\gamma$}
We use Theorem~\ref{thm:purpernonunit}   to deduce several conditions
for pure periodicity in ${\mathbb Q}(\beta)$. 
That 0 is an inner point of the complete central tile $\Tt$ yields
a first sufficient condition for a rational number to have purely periodic
expansion. We can see this property as a  generalisation of
Corollary~\ref{gamabetaunit}.  

\begin{corollary}
Let $\beta$ be a Pisot  number  that satisfies the  finiteness property
(F). There exist $m$ and $v$ such that for every   $x = \frac{N(\beta)^m p}{q} \in
{\mathbb Q}$, with $\gcd(N(\beta),q)=1$,  and $x \leq v$, then
$x\in\Pi_\beta^{(r)}$.   
\end{corollary}

\begin{proof}
Let $M$ be the maximum of the $\vert N(\beta)\vert_v$, for
$v\in\Sb$, $v$ non-Archimedean. We have $M<1$. Therefore, for
$x=\frac{N(\beta)^m p}{q} \in {\mathbb Q}$, with $\gcd(N(\beta),q)=1$ and
$x\leq v$, we have $\Vert x\Vert\leq \max(v, M^m)$. Since 0
is an inner point of $\Tt$ - actually an inner point of $\Ta{1}$ by
(\ref{eq:lesboules}) - it follows that $(-\phb(x),x)\in\Ta{1}\times[0,1)$
if $m$ is big enough, and $v$ small enough.
\end{proof}

\subsection{From  the topology of the central  tile to that of
 $\Pi_\beta^{(r)}$}  

We begin with completing the notation introduced in
Sections~\ref{subsec:bourbaki} and~\ref{subsec:crs}. 
\begin{notation}
Recall that $(\beta)=\prod_{i=1}^\nu\mathfrak{P}_i^{n_i}$, and that $\GK_f$ is
the product of the corresponding local fields. We denote by $\phi_f$ the
associated embedding, so that 
 $$\phi_\beta(x)=(\phi_\infty(x),\phi_f(x))\in\GK_\infty\times\GK_f.$$
  We also write $\OO_f=\prod_{i=1}^\nu\OO_{\mathfrak{P}_i}$ and  we denote by
  $\OO_{(\beta)}$ its reciprocal image by $\phi_f$, that is, 
    $$\OO_{(\beta)}=\{x\in\GQ(\beta);\, \forall i,\, 1\leq i\leq  \nu:\,
    v_{\mathfrak{P}_i}(x)\geq 0\}\supset \OO.$$ There are primes $p_i$ such
    that $\mathfrak{P }_i\cap\GZ=p_i\GZ$. We write
    $\mathbf{k}_f=\prod_{i=1}^\nu\GQ_{p_i} $ and
    $\mathfrak{o}_f=\prod_{i=1}^\nu\GZ_{p_i}$. We also introduce
\[
\GZ_{(N(\beta))}=\left\{\frac pq;\,  \gcd(q,N(\beta))=1\right\}.
\]
If $N(\beta)$ is  prime, this notation coincides with the usual one concerned
with localisation. Let us finally introduce the canonical projections 
\[
\pi_\infty\colon\GK_\beta\rightarrow\GK_\infty
\mbox{ and } \pi_f\colon\GK_\beta\rightarrow\GK_f.
\] 
\end{notation}
Most of this notation is summarised in the commutative
diagrams below: 
\begin{equation}\label{eq:commutativediagram}
\begin{CD}
\GQ(\beta) @>\phi_f >>\displaystyle{\prod_{i=1}^\nu\GK_{\mathfrak{P}_i}}=\GK_f\\
@A i_1 AA @AA i_2 A\\
\GQ @> \phi_f >>\displaystyle{\prod_{i=1}^\nu\GQ_{p_i}=\mathbf{k}_f}
\end{CD}\qquad\qquad
\begin{CD}
\OO_{(\beta)} @>\phi_f >>\displaystyle{\prod_{i=1}^\nu\OO_{\mathfrak{P}_i}=\OO_f}\\
@A i_1 AA @AA i_2 A\\
\GZ_{(N(\beta))} @> \phi_f
>>\displaystyle{\prod_{i=1}^\nu\GZ_{p_i}=\mathfrak{o}_f} 
\end{CD}
\end{equation}

 We want to generalise the idea of the proof of Corollary  \ref{gamabetaunit}.
 Since we are from now one interested in the beta-expansion of rational
 integers, our first goal is to understand how they imbed into $\GK_\beta$. The
 Archimedean embedding is trivial: $\phi_{\infty}(x)=(x,x,\ldots,x)$. 

\begin{notation}[Diagonal sets] Let $A\subset\GR$. Then
  $$\Delta_{\infty}(A):=\{(x,\ldots,x);\, x\in A\}\subset\GK_\infty$$   stands
  for the    set of $(r+s-1)$-uples of elements of $A$ whose coordinates are
  all equal. 

By an abuse of language, when $A$ is reduced to a single point $A = \{a \}$, we
will use the notation $\Delta_{\infty}(a)$ for the point
$\phi_\infty(a,\ldots,a)\in\GK_\infty$.   
\end{notation}

We now need to understand   the  non-Archimedean   embedding  of
$\GZ_{(N(\beta))}$: this is  the object of  
 Lemma~\ref{lem:localsurjective} and Proposition~\ref{prop:charactdensity}
 below. 

\begin{lemma}\label{lem:localsurjective}
Let $V$ be a non-empty open subset of $\GZ_{(N(\beta))} $. Then
\[
\overline{\phi_f(V)}=\overline{\phi_f(\GZ_{(N(\beta))}
)}\quad\textrm{and}\quad
\overline{\phi_\beta(V)}=\Delta_{\infty}(\overline{V}) \times
\overline{\phi_f(\GZ_{(N(\beta))})}.
\]
Furthermore, for any non-empty interval $I$ in $[0,1]$, we have
\begin{equation}\label{eq:phiofi}
\overline{\phb(I \cap \GZ_{(N(\beta))})}=\Delta_{\infty}(\bar
I)\times\overline{\phi_f(\GZ_{(N(\beta))})}.
\end{equation}
The same results hold  if one replaces $\GZ_{(N(\beta))} $ by $\GQ$,
 $\GQ(\beta) $, or $\OO_{(\beta)} $.
\end{lemma}

\begin{proof} We only prove the result for $\GZ_{(N(\beta))} $ - the other
  cases being similar. Let $V$ be a non-empty open subset of
$\GZ_{(N(\beta))}
  $ and $u\in V$. For $y\in \overline{\phi_f(\GZ_{(N(\beta))})}$, there
exists
  a sequence $(x_n)_n$ in $\GZ_{(N(\beta))} $ such that
  $\lim\phi_f(x_n)=y-\phi_f(u)$ (using that $\phi_f$ is an additive group
  homomorphism). Let us introduce $\vartheta_n=(1+N(\beta)^n)^{-1}$. Then
  $\vartheta_n\in\GZ_{(N(\beta))} $, and we have both
$\lim\vartheta_n=0$ and
  $\lim\phi_f(\vartheta_n)=1$. Then we can choose a subsequence
$(\sigma(n))_n$
  such that $u+x_n\vartheta_{\sigma(n)}\in V$ and $\lim
  \phi_f(u+x_n\vartheta_{\sigma(n)})=y$. Finally, $\phi_{\beta}(x_n)$
converges to $(\phi_{\infty}(u),y)$.  

This means that $\Delta_{\infty}({V}) \times
\overline{\phi_f(\GZ_{(N(\beta))}
  )} \subset \overline{\phi_\beta(V)}$; taking the closure yields
$\Delta_{\infty}(\overline{V}) \times \overline{\phi_f(\GZ_{(N(\beta))} )}
\subset \overline{\phi_\beta(V)}$. We conclude by noticing that the
definition
of $\phi_\beta$ directly ensures that $\overline{\phi_\beta(V)} \subset
\Delta_{\infty}(\overline{V}) \times
\overline{\phi_f(\GZ_{(N(\beta))})}$. Hence
we have proved
$\overline{\phi_\beta(V)}=\Delta_{\infty}(\overline{V})
\times\overline{\phi_f(\GZ_{(N(\beta))})}$. The equality
  $\overline{\phi_f(V)}=\overline{\phi_f(\GZ_{(N(\beta))})}$ follows by
applying the  projection $\pi_f$.

Equation~(\ref{eq:phiofi}) is clearly satisfied if $I$ is open. In
general, it follows from
\[
\Delta_{\infty}(\bar I)\times\overline{\phi_f(\GZ_{(N(\beta))})} =
\overline{\phb(\mathring{I} \cap \GZ_{(N(\beta))})} \subset\overline{\phb(I
  \cap \GZ_{(N(\beta))})} \subset\Delta_{\infty}(\bar
I)\times\overline{\phi_f(\GZ_{(N(\beta))})}.
\]
\end{proof}
We deduce from Lemma \ref{lem:localsurjective} and
Theorem~\ref{thm:purpernonunit}:  

\begin{corollary} \label{cor:stripe}
Let  $0<\varepsilon \leq\min\{T_{\beta}^{a-1}(1), \, a \in \{1, \dots, n \}\}$.
Then  $\gamma(\beta) \geq \varepsilon$ if  
$\Delta_{\infty} ([0,\varepsilon]) \times \overline{\phi_f(\GZ_{(N(\beta))})}
\subset  -\widetilde{\mathcal{T}}$. 
\end{corollary}

Let us stress the fact that there is no reason here for   $\phb(\GQ)$ 
to be dense  in $\GK_\beta$, contrarily to  what happens  for the Archimedean
part. 

\begin{proposition}\label{prop:charactdensity}
The following propositions are equivalent:
\begin{enumerate}
\item Let $0\leq u<v<\infty$. Then
$\overline{\phi_f (\GZ_{(N(\beta))}\cap(u,v))}=\OO_f$;
\item The set $\phi_f(\GQ)$ is dense in $\GK_f$;
\item For all $i$, $1\leq i\leq \nu$, we have
  $e(\mathfrak{P}_i)=f(\mathfrak{P}_i)=1$, and the prime numbers $p_i$ are all
  distinct. 
\item The norm $N(\beta)$ is square-free and none of its prime divisors
  ramifies. 
\end{enumerate}
\end{proposition}
\begin{proof}  For given $i$, one has  
  $[\GK_{\mathfrak{P}_i}:\GQ_{p_i}]=e(\mathfrak{P}_i)f(\mathfrak{P}_i)$. By
  completeness of the $p$-adic fields (resp. $p$-adic rings) the image by 
$i_2$ of $\mathbf{k}_f$ (resp. $\mathfrak{o}_f$) in the commutative
diagram~(\ref{eq:commutativediagram}) is closed in $\GK_f$
(resp. $\OO_f$). Hence, these images are dense if and only if those products
are equal, that is if $\GK_{\mathfrak{P}_i}=\GQ_{p_i}$ for all $i$, i.e.,
$e(\mathfrak{P}_i)=f(\mathfrak{P}_i)=1$.

Moreover, by the Chinese remainder
theorem, the image by $\phi_f$ of $\GQ$ (resp. $\GZ_{(N(\beta))}$) is dense in
$\mathbf{k}_f$ (resp. $\mathfrak{o}$) if and only if
the $p_i$ are distincts. Hence we have proved that (2), as (1) for
$\kappa=+\infty$, are equivalent to (3). The equivalence with (1) with an
arbitrary non-empty open interval $(u,v)$ is given by
Lemma~\ref{lem:localsurjective}.  

Finally, the equivalence of (3) and (4) follows from
\[
N(\beta)=N((\beta))=\prod_{i=1}^\nu N(\mathfrak{P}_i)=\prod_{i=1}^\nu 
p_i^{f(\mathfrak{P}_i)}.
\]
\end{proof} 

\begin{remark} If the prime numbers $p_i$ are not distinct, there is a
  partition of $\nu$, $\nu=\mu_1+\cdots +\mu_\ell$ and a suitable reordering of
  the prime ideals $\mathfrak{P}_1,\ldots, \mathfrak{P}_\nu$ containing
  $\beta$,  such that one has the equality of multisets 
\[
\{p_1,p_2,\ldots, p_\nu\}=\{\underbrace{p_1,\ldots,p_1}_{\# \mu_1},
\underbrace{p_2,\ldots,p_2}_{\# \mu_2},\ldots,
\underbrace{p_\ell,\ldots,p_\ell}_{\# \mu_\ell}\}. 
\] 
Then, $\overline{\phi_f(\GQ)}$ (resp. $\overline{\phi_f(\GZ_{(N(\beta))})}$) is
equal to $\prod_{j=1}^\ell\Delta(\GQ_{p_j}^{\mu_j})$ (resp.
$\prod_{j=1}^\ell\Delta(\GZ_{p_j}^{\mu_j})$), where $\Delta(M^\mu)$ denotes the
set of $\mu$-uples of elements of $M$ whose coordinates are all equal. 
\end{remark}

\subsection{Topological properties of  $\Pi_\beta^{(r)}$}
Before  being able to  deduce  bounds  on $\gamma(\beta)$ from
Corollary~\ref{cor:stripe}, we need to preliminary investigate the topological
structure of $\Pi_\beta^{(r)}$.

We already know that $\Pi_\beta^{(r)}\subset\GZ_{(N(\beta))}$ by Lemma
\ref{lem:pp}. We endow $\Pi_\beta^{(r)}$ with the induced topology of $\GR$ on
$\GZ_{(N(\beta))}$. The following proposition investigates the extremities of
$\Pi_\beta$'s connected components. An example of such a component is of course
$[0,\gamma(\beta)]$ (or $[0,\gamma(\beta))$).

\begin{theorem}\label{theo:uv}
Let $(u,v)$ be a non-empty open interval with
$(u,v)\cap\GZ_{(N(\beta))}\subset\Pi_\beta^{(r)}$. 

If $v\in\GZ_{(N(\beta))}$, then $v\in\Pi_\beta^{(r)}$.  
If the assumptions of Proposition~\ref{prop:charactdensity} are satisfied  and $v\in {\mathbb Q}$,
then the same conclusion $v\in\Pi_\beta^{(r)}$ holds.

If $(u,v)$ as above is maximal and $v<1$, then there are three possibilities
for $v$, namely:  
\begin{enumerate}
\item[(A)] There exists $a\in\A$ such that
\[
\Delta_\infty(v)\in\pi_\infty
\left(-\Ta{a}\cap\left(\Delta_\infty(T_\beta^{(a-1)}(1))
  \times\overline{\phi_f(\GZ_{(N(\beta))})}\right)\right). 
\]
In particular, $v=T_\beta^{(a-1)}(1)$.
\item[(B)] There exist $a$ and $b$ in $\A$ such that 
\[
\Delta_\infty(v)\in\pi_\infty\left(-\Ta{a}\cap -\Ta{b}\cap
  \left(\Delta_\infty([T_\beta^{(b-1)}(1), T_\beta^{(a-1)}(1)))\times
  \overline{\phi_f(\GZ_{(N(\beta))})}\right)\right).  
\]
In particular, $T_\beta^{(b-1)}(1)\leq v<T_\beta^{(a-1)}(1)$.
\item[(C)] There exist $a\in\A$ and $x\in\GZ[1/\beta]\cap(0,1)$ such that
\[
\Delta_\infty(v)\in\pi_\infty\left(-\Ta{a}\cap -\Tt(x)\cap
  \left(\Delta_\infty((0, T_\beta^{(a-1)}(1)))\times
    \overline{\phi_f(\GZ_{(N(\beta))})} \right)\right).
\]
In particular, $v<T_\beta^{(a-1)}(1)$.
\end{enumerate}
Cases (B) and (C) are not exclusive of each other.  The same results hold for
$u$, $u>0$.
\end{theorem}
\begin{proof} Let $(u,v)$ be a non-empty open interval with
$(u,v)\cap\GZ_{(N(\beta))}\subset\Pi_\beta^{(r)}$. Assume that
$v\in\GZ_{(N(\beta))}$. Then, by Lemma~\ref{lem:localsurjective}, one
can construct a sequence $(z_n)_n$ in $(u,v)$ such that $\lim z_n=v$ and
$\lim\phi_f(z_n)=\phi_f(v)$. Furthermore,
$\lim z_n=v$ is equivalent to $\lim\phi_\infty(z_n)=\phi_\infty(v)$.
Hence, we have $\lim\phb(z_n)=\phb(v)$. Moreover, by taking a subsequence, we
may assume that for some $a\in\A$, one has $(-\phb(z_n), z_n)\in\Ta{a}\times[0,
T_\beta^{(a-1)}(1))$ for all $n$. Then $(-\phb(v), v)\in\Ta{a}\times[0,
T_\beta^{(a-1)}(1)]$. By Lemma~\ref{lem:Torbit}, the assumption
$v\in\GZ_{(N(\beta))}\subset\GQ$ guarantees that $v\neq
T_\beta^{(a-1)}(1)$. Therefore, we have that $(-\phb(v), v)\in\Ta{a}\times[0,
T_\beta^{(a-1)}(1))$ and $v\in\Pi_\beta$. The same argument applies to
$v\in\GQ$ under the assumptions of Proposition~\ref{prop:charactdensity}. The
case of $u$ is similar.\\

We now assume that the interval $(u,v)$ is maximal and $v\neq 1$. We first
claim that there 
exists a sequence $(y_n)_n$ in $\GZ_{(N(\beta))}\setminus\Pi_\beta^{(r)}$ with
$\lim y_n=v$. By the maximality of $(u,v)$, it is trivial, but if there exists
$w>v$ such that $(u,w)\cap\GZ_{(N(\beta))}\subset\Pi_\beta^{(r)}$ and
$v\in\GZ_{(N(\beta))}\setminus\Pi_\beta^{(r)}$. By the first part of the
theorem, this cannot happen, and our claim is proved.

Let us then start with a sequence $(y_n)_n$ with
$v<y_n$, $\lim y_n=v$ and $y_n\not\in\Pi_\beta$. By compacity, one may
assume that $(\phi_\beta(y_n))_n$ converges, to $(\Delta_\infty(v), z)$,
say, with $z\in \overline{\phi_f(\GZ_{(N(\beta))})})$. By
Lemma~\ref{lem:localsurjective}, there exists a sequence $(z_n)_n$ with
$u<z_n<v$, $\lim z_n=v$ and $\lim\phb(z_n)=(\Delta_\infty(v), z)$. By
extracting a subsequence, we also may assume that there exists $a\in\A$ with
$(\phb(z_n),z_n)\in-\Ta{a}\times[0,T_\beta^{(a-1)}(1))$ for all $n$. The first
possibility to take in account is $v=T_\beta^{(a-1)}(1)$. Since $\Ta{a}$ is
closed, we then have $(\Delta_\infty(v),z,v)\in
-\Ta{a}\times\{T_\beta^{(a-1)}(1)\}$. In other words, one gets the possibility
(A) of the theorem:
\[
\Delta_\infty(v)\in\pi_\infty
\left(-\Ta{a}\cap\left(\Delta_\infty(T_\beta^{(a-1)}(1))
  \times\overline{\phi_f(\GZ_{(N(\beta))})}\right)\right). 
\]

From now on, we may assume that $v\neq T_\beta^{(a-1)}(1)$ (that does not mean
that $v$ could not be equal to an other element of the $T_\beta$-orbit of
1). We then have $(\Delta_\infty(v), y)\in
-\Ta{a}\times[0,T_\beta^{(a-1)}(1))$. Since 
$y_n\not\in\Pi_\beta$, we get $\phb(y_n)\not\in-\Ta{a}$. For fixed $n$, there
are two possibilities:
\begin{enumerate}
\item[(i)]  $\phb(z_n)\in -\Tt$. Since $z_n\not\in\Pi_\beta$, we have that
  $\phb(z_n)\in   -\Ta{b}$ for some $b\in\A$ such that $T_\beta^{(b-1)}(1)\leq
  v$. 
\item[(ii)] $\phb(z_n)\not\in-\Tt$. Then, by Proposition~\ref{prop:covering},
  there  exists $x_n\in\GZ[1/\beta] \cap(0,1)$ such that $\phb(z_n) \in -\Tt
  (x_n)$. 
\end{enumerate}
At least one of the properties (i) of (ii) has to be satisfied for infinitely
many $n$'s.

If that is the case for (i), since $\A$ is finite, there is a $b$
corresponding to a further subsequence of $(z_n)_n$. Taking the limit, we get
$\phb(v)\in -\Ta{b}$. Hence case (B) of the theorem:
\[
\Delta_\infty(v)\in\pi_\infty\left(-\Ta{a}\cap -\Ta{b}\cap
  \left(\Delta_\infty([T_\beta^{(b-1)}(1), T_\beta^{(a-1)}(1)))\times
  \overline{\phi_f(\GZ_{(N(\beta))})}\right)\right).  
\]

If there are infinitely many $n$'s satisfying (ii),
Proposition~\ref{prop:covering} shows that the family
$\{x_n,n\in\GN\}$ is finite. Hence, by extracting a subsequence, there
is some $x\neq 0$ with $\phb(z_n) \in -\Tt (x)$. Taking the limit, we get case
(C): 
\[
\Delta_\infty(v)\in\pi_\infty\left(-\Ta{a}\cap -\Tt(x)\cap
  \left(\Delta_\infty((0, T_\beta^{(a-1)}(1)))\times
    \overline{\phi_f(\GZ_{(N(\beta))})} \right)\right).
\]
\end{proof}

\begin{proposition}
If the finitess  property (F) is satisfied, then $\Pi_\beta$ is dense in
$\GZ_{(N(\beta))}$.
\end{proposition}
\begin{proof}
If the property (F) is satisfied, then $\Ta{1}$ contains a neighbourhood
of the origin, hence $\overline{\phi_f(\Int(\beta))}$ contains
$N(\beta)^m\overline{\phi_f(\GZ_{(N(\beta))})}$ for some $m\geq 0$. Then
\[
\phb^{-1}\left(\Delta_\infty([0,1))\times
N(\beta)^m\overline{\phi_f(\GZ_{(N(\beta))})}\right)\subset \Pi_\beta,
\]
and is dense by Lemma~\ref{lem:localsurjective}.
\end{proof}

\subsection{Upper  and lower bounds for  $\gamma(\beta)$}
We now  have collected all the required material to be able to deduce
 upper and lower bounds for  $\gamma(\beta)$. The present section collects
 results that may be of some interest in   every dimension, whereas
 Section~\ref{subsec:quad} is devoted to the quadratic case.\\

A first  upper bound for $\gamma(\beta)$ can be directly deduced from
Theorem~\ref{thm:purpernonunit}. We consider the intersection between the
complete central subtiles and the set of points whose canonical Archimedean
projection by $\pi_{\infty}$ belong to the  diagonal  sets  of the form 
$\Delta_{\infty}([0, T_{\beta}^{(a-1)}(1)))$.
\begin{proposition}\label{lem:firstupperbound}
 Let  $\beta$ be a Pisot number.  One has:
\[  
\gamma(\beta) \leq 
\max \left  \{ T_{\beta}^{(a-1)}(1); \   a\in\A, \
  (-\Ta{a})\cap\pi_\infty^{-1}\Delta_{\infty}([0,
  T_{\beta}^{(a-1)}(1)))\neq\emptyset  \right\}.
\]
\end{proposition}

\begin{proof}
Let $x \in\GQ\cap[0,1)$. If  $(\phb(x),x)$  belongs to $\bigcup_{a \in \A}
(-\Ta{a})\times[0, T_{\beta}^{a-1}(1))$,  then   there exists $a \in \A$ such
that $ \pi_{\infty} \circ\phb(x) \in -\Ta{a}\cap  \Delta_{\infty}([0,
    T_{\beta}^{(a-1)}(1)).$  Hence if
\[
x > \max \left\{
  T_{\beta}^{a-1}(1); \,a \in \A, \, (-\Ta{a})\cap\pi_\infty^{-1}
  \Delta_{\infty}([0, 
    T_{\beta}^{(a-1)}(1)))\neq \emptyset\right\},
\]
then $(\phb(x),x)$ does not belong to $\bigcup_{a \in \A} (-\Ta{a})\times[0,
T_{\beta}^{a-1}(1)).$  We deduce from Theorem~\ref{thm:purpernonunit}  that
its $\beta$-expansion is not purely periodic.   
\end{proof}

Let us stress the point that this upper bound is quite rough: if the finiteness
property (F) is satisfied, then the inequality yields the trivial bound
$\gamma(\beta) \leq 1$. Indeed Proposition~\ref{prop:innerpoint} 
says that $\Ta{1}$ contains a neighbourhood of the origin. Hence the
intersection $(-\Ta{1})\cap\pi_\infty^{-1} \Delta_{\infty}([0, T_{\beta}(1)))$
is not empty, which yields $\gamma(\beta)\leq 1$.\\

However, Theorem~\ref{thm:purpernonunit} states that real numbers have a purely
periodic expansion if their embedding is included in the representation
$\bigcup_{a \in \A} (-\Ta{a})\times[0, T_{\beta}^{a-1}(1))$ of the natural
extension of $T_{\beta}$. From Lemma~\ref{lem:localsurjective}, we know that an
interval of rationals $(\eta,\nu)\cap \GZ_{(N(\beta))}$ embeds in $\GK_{\beta}$
as the product of a diagonal set with a local part whose closure is independant
of $(\eta,\nu)$. We deduce below a recursive characterisation for
$\gamma(\beta)$.     

\begin{notation}\label{not:order}
Let us order  and relabel  the elements in $\A$
as follows:  we  set 
$\A=\{a_1,\ldots,a_n\}$   with
\[
T_\beta^{a_1-1}(1)<T_\beta^{a_2-1}(1)<\cdots
<T_\beta^{a_{n-1}-1}(1)<T_\beta^{a_n-1}(1)=1.  
\]
Clearly, $a_n=1$. For notational convenience, we state $T_\beta^{a_0-1}(1)=0$. 
\end{notation}

\begin{proposition} \label{prop:rec}
Let $\beta$ be a Pisot number.
\begin{itemize}
\item $\gamma(\beta)\geq T_\beta^{a_{k}-1}(1)$ if and only if:
\[ 
\gamma(\beta)\geq T_\beta^{a_{k-1}-1}(1)\quad\textrm{and}\quad\Delta_{\infty}([
  T_\beta^{a_{k-1}-1}(1),T_\beta^{a_{k}-1}(1)])\times
  \overline{\phi_f(\GZ_{(N(\beta))})}\subset\bigcup_{j=k}^n(-\Ta{a_j}).
\]
\item If $T_\beta^{a_{k-1}-1}(1)< \gamma(\beta)\leq T_\beta^{a_{k}-1}(1)$, then 
\begin{equation}\label{eq:explicitgamma1}
\gamma(\beta)=\sup\left\{\eta\geq T_\beta^{a_{k-1}-1}(1);\; \Delta_{\infty}([
  T_\beta^{a_{k-1}-1}(1),\eta])\times
  \overline{\phi_f(\GZ_{(N(\beta))})}\subset\bigcup_{j=k}^n(-\Ta{a_j})\right\}.
\end{equation}
\end{itemize}
In particular, if $\Tt$ does not contain
$\Delta_{\infty}([0,\eta])\times\overline{\phi_f(\GZ_{(N(\beta))})}$ for any
positive $\eta$, then $\gamma(\beta)=0$. 
\end{proposition} 
\begin{proof} Let $I$  a non-empty open interval  in $[0,1]$.
By Lemma~\ref{lem:localsurjective}, $I\cap \GZ_{(N(\beta))}\subset\Pi_\beta$
if and only if 
\[
\Delta_{\infty}({\bar
  I})\times\overline{\phi_f(\GZ_{(N(\beta))})}
\subset\bigcup_{j=k}^n(-\Ta{a_j}).
\]
Equation~(\ref{eq:explicitgamma1}) follows from~(\ref{eq:phiofi}) and
Theorem~\ref{thm:purpernonunit} too. The last assertion is a particular case of
~(\ref{eq:explicitgamma1}) when $k=1$ and of the observation that
$\phi_f(\GZ_{(N(\beta))})=-\phi_f(\GZ_{(N(\beta))})$. 
\end{proof}
%\comment{We know that $\phi_f(\OO)$ is dense in
%  $\prod_{i=1}^\nu\OO_{\mathfrak{P}_i}$. Does it hold for
%  $\phi_f(\Int(\beta))$? Or at least for $\phi_f(\GZ[\beta])$? If not, then
%  $\gamma(\beta)=0$ (under the assumptions of
%  Proposition~\ref{prop:charactdensity}). We need to find an example of
%  application. {\bf A priori je n'ai pas trouve d'exemple utilisable pour
%  l'instant. en dimension 
%2 je pense qu'il n'y a pas d'exemples de ce type. Il faudrait chercher en
%degre 3 ou 4 pour trouver quelque chose, demander 
%a Sing peut etre}}  

This result has a geometric interpretation related to the natural extension of
$T_{\beta}$. Denote by $\Delta$ the diagonal line in ${\mathbb
  K}_{\infty}\times{\mathbb R}$, that is, the Euclidean component of the
natural extension. Proposition \ref{prop:rec} means  that $\gamma(\beta)$ is
the largest part of  $\Delta$ starting from $0$ such that its product with the
full non-Archimedean component $\overline{\phi_f(\GZ_{(N(\beta))})}$ is totally
included in the  natural extension $\bigcup_{a \in \A} (-\Ta{a})\times[0,
T_{\beta}^{a-1}(1))$.   

\begin{figure}[ht]
\begin{center}
\includegraphics[height=5cm]{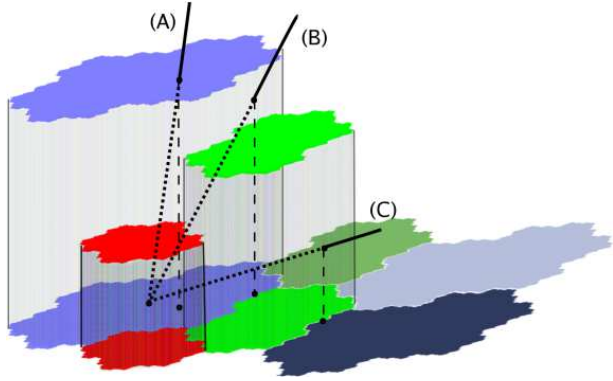}
\caption{ Illustration of the three cases of  Theorem \ref{theo:uv} and
  Proposition \ref{prop:rec}. 
We have chosen a unit Pisot number for the illustration  of these
three cases for the sake of  clarity.  By Proposition \ref{prop:rec},
$\gamma(\beta)$ is given by the largest  part of the diagonal line to be
fully included in the natural extension 
$\bigcup_{a \in \A} (-\Ta{a})\times[0, T_{\beta}^{(a-1)}(1))$. The natural
extension is represented with subtiles $-\Ta{a}$ in the horizontal direction,
and the interval $[0, 1)$ in the vertical axis. 
Then, the natural extension consists in union of cylinders with fractal
horizontal base and vertical height. The height of 
the cylinder with basis  $-\Ta{a}$ is  $T_{\beta}^{(a-1)}(1)$. Depending on the
location of the point where the diagonal first goes out from the natural
extension, we get the different situations unearthed in Theorem~\ref{theo:uv}. 
\newline  
Situation (A) means that $\gamma(\beta)$  belongs  to the orbit of $1$ 
under the action of $T_{\beta}$ and that its Euclidean embedding
$\phi_{\infty}(\gamma(\beta))$ simultaneously
is the Euclidean part of a point of the corresponding subtile. 
Then, the diagonal starts from 0 and exits from the natural extension 
on a plateau with height $T^{(a-1)}_{\beta}(1)$. 
\newline
Situation (B)  involves the   intersection between two  complete central
subtiles tiles $(-\Ta{a})\cap(-\Ta{b})$.  The diagonal line  goes out from the
natural extension on a vertical line above the intersection 
between two   subtiles. The main point is that the plateau of the lowest
cylinder ($T^{(b-1)}(1)$) lies below the diagonal line  whereas the plateau of
the  upper cylinder ($T^{(a-1)}(1)$) lies above it.
\newline
Situation (C) means that the diagonal line completely crosses the natural
extension and exits above a new $x$-tile.  
}\label{Fig:diago} 
\end{center}
\end{figure}

In the unit case, since the representation contains only Archimedean
components, Proposition \ref{prop:rec} simply means 
that $\gamma(\beta)$ is the length of the largest diagonal interval that is
fully included in the natural extension (see an illustration in
Fig. \ref{Fig:diago}).

Theorem~\ref{theo:uv} offhands yields  lower and upper bounds for $\gamma(\beta)$.
\begin{proposition}\label{prop:uvbounds}
We introduce some local notation. For $a$ and $b$ in $\A$ such that
$T_\beta^{(b-1)}(1)\leq T_\beta^{(a-1)}(1)$, let
\[
A_{a,b}=\pi_\infty\left(-\Ta{a}\cap -\Ta{b}\cap
  \left(\Delta_\infty([T_\beta^{(b-1)}(1), T_\beta^{(a-1)}(1)])\times
  \overline{\phi_f(\GZ_{(N(\beta))})}\right)\right)\subset\GK_\infty.  
\]
For $a\in\A$ and $x\in\GZ[1/\beta]$, let
\[
B_{a,x}=\pi_\infty\left(-\Ta{a}\cap -\Tt(x)\cap
  \left(\Delta_\infty((0, T_\beta^{(a-1)}(1)))\times
    \overline{\phi_f(\GZ_{(N(\beta))})} \right)\right)\subset \GK_\infty.
\]
Finally, let 
\[
A=\bigcup_{\genfrac{}{}{0pt}{1}{(a,b)\in\A^2}{T_\beta^{(b-1)}(1)\leq   
      T_\beta^{(a-1)}(1)}}A_{a,b}\quad\textrm{and}\quad
B=\bigcup_{\genfrac{}{}{0pt}{1}{a\in\A}{x\in\GZ[1/\beta]\cap(0,1)}}B_{a,x}
\]
Then, an upper bound for  $\gamma(\beta)$ is given by
\begin{equation}\label{eq:minorationforgamma}
\gamma(\beta)\geq
 \min\left(\min_{\substack{(a,b)\in\A^2\\ T_\beta^{(b-1)}(1)\leq
      T_\beta^{(a-1)}(1)\\ A_{a,b}\neq\emptyset}}\min_{x\in A_{a,b}}\Vert
 \pi_\infty(x)\Vert_\infty,\min_{\substack {a\in\A\\
     x\in\GZ[1/\beta]\cap(0,1)\\ 
    B_{a,x}\neq\emptyset}}\inf_{x\in B_{a,x}}\Vert\pi_\infty(x)\Vert_\infty \right).
\end{equation}
 A lower bound for  $\gamma(\beta)$ is the following:
\begin{equation}\label{eq:majorationforgamma}
\gamma(\beta)\leq\max\left\{\eta;\,
  [0,\eta]\subset A\cup B\right\}. 
\end{equation}
\end{proposition}

\begin{proof}
First note that the infimum in \eqref{eq:minorationforgamma} is due to the fact
that $B_{a,x}$ does not need to be compact. We use Theorem~\ref{theo:uv} and
the fact, that, by definition, $\gamma(\beta)$ is the largest number $\gimel$
such that $(0, \gimel) \cap\GZ_{(N(\beta))}\subset\Pi_\beta^{(r)}$.  
Situation (A) in Theorem~\ref{theo:uv} implies that there exists $a \in
{\mathcal A}$ such that $\gamma(\beta) = T_{\beta}^{(a-1)}(1)$ and
$\gamma(\beta)\in\pi_\infty (-\Ta{a})$; it reads off that $\gamma(\beta)\in
A_{a,a}$.  Situation (B) implies that there exist $a,b \in {\mathcal A}$ with
$T_{\beta}^{b-1}(1)<T_{\beta}^{a-1}(a)$ such that $\gamma(\beta) \in
A_{a,b}$. However, the interval is closed in the present proposition, as it is
half-closed in Theorem~\ref{theo:uv}. Nevertheless, by continuity of
$\pi_\infty$, taking open or closed intervals in $B_{a,x}$ or $A_{a,b}$ has no
influence on the infimum we are interested in. Situation (C) reads off that
there exist $a \in {\mathcal A}$ and $x \in\GZ[1/\beta]\cap(0,1)$ such that
$\gamma(\beta) \in B_{a,x}$. Since one the 3 situations must occur, we deduce
that $\gamma(\beta)$ is greater than the smallest of the infimum of all these
sets. Formulas~\eqref{eq:majorationforgamma}  hold for the same reasons. 

The three cases are illustrated in Fig.~\ref{Fig:diago}.
\end{proof}

\subsection{Quadratic Pisot numbers}\label{subsec:quad}
Let us now consider  the particular case of quadratic Pisot numbers of
degree 2, for which many things can be done explicitely. For instance,
$\GQ(\beta)$ is an extension of degree two, and then the
  algebraic conditions (3) or (4) of Proposition~\ref{prop:charactdensity} can
  be easily tested. Indeed, let $d$ be  the square-free positive rational
  integer such that   $\GQ(\beta)=\GQ(\sqrt d)$. Then the discriminant
  $\delta_{\GQ(\beta)}$ of the quadratic field is $d$ if $d\equiv 1\pmod
4$ and   $4d$ if $d\equiv 2,3 \pmod 4$.
\begin{corollary}\label{cor:quadraticdensity}
 If $\beta^2= a\beta+b$, with
  $(a,b)\in\GZ^2$, $b\neq 0$, the  equivalent conditions of
  Proposition~\ref{prop:charactdensity} are satisfied if and only if:
\begin{enumerate}
\item $b$ is square free,
\item $b$ is coprime with $\delta_{\GQ(\beta)}$,
\item $d$ is a  quadratic residue with respect to all odd prime divisors of
$b$,
\item $d\equiv 1\pmod 8$ if $b$ is even.
\end{enumerate}
\end{corollary}

The Euclidean representation space ${\mathbb K}_{\infty}$
is a one-dimensional line. Consequently, the diagonal
$\Delta_{\infty}([0,\varepsilon])$ is indeed the interval $[0,\varepsilon]
\subset\GK_{\infty}=\GR$. This allows  us to use graphical representation of
the complete central tile to conjecture lower bounds  for $\gamma(\beta)$.\\

A particularly manageable case is the following: $(\beta)=\beta\OO$ is a
prime ideal lying above a prime number $p$, that splits. Hence $(\beta)$ has
inertia degree 1, we have $N((\beta))=\vert N(\beta)\vert=p$, and
$\GK_\beta=\GR\times\GQ_p$ (that is a special case of
Corolary~\ref{cor:quadraticdensity}). We can represent $\GZ_p$ by the Mona
map $\GZ_p\ni x = \sum a_i p^i\mapsto \sum a_i p^{-i} \in [0,1]$. This
mapping
is onto, continuous and preserves the Haar measure, but is it not a morphism
for the addition. Corollary~\ref{cor:stripe} implies that $\gamma(\beta)\geq
\varepsilon$ if and only if a stripe of length
$\varepsilon$ is totally included in the representation of the central
tile, as
illustrated by Fig. \ref{fig:10_3} and Fig. \ref{fig:4_3}  below.

\begin{figure}
\includegraphics[height=3cm, width=6cm]{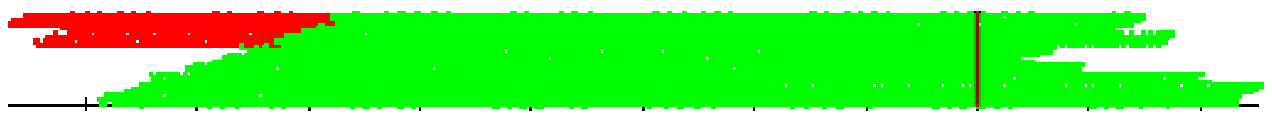}
\includegraphics[height=2.5cm]{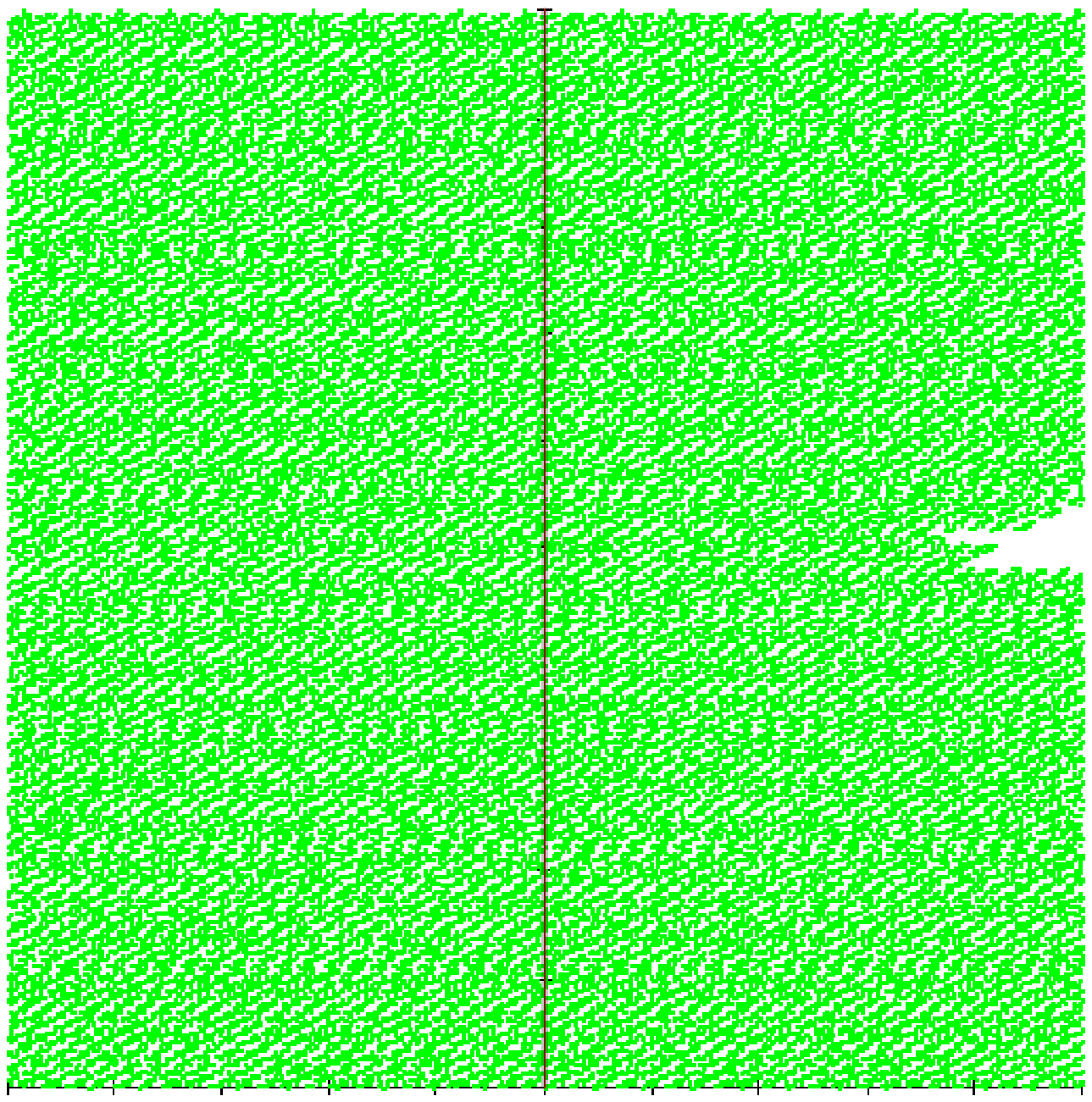}
\caption{A representation of the complete central tile for $\beta=5+2\sqrt
  7$. Then $\beta$ has minimum polynomial $X^2-10X-3$ and
$N((\beta))=-3$. The
  quadratic field $\GQ(\beta)=\GQ(\sqrt 7)$ has discriminant 28, hence
  $\left(\frac{28}{3}\right)=\left(\frac 13\right)=1$: the discriminant is a
  quadratic residue modulo 3 and Corollary~\ref{cor:quadraticdensity} shows
  that the complete  central tile is a subset of
${\mathbb R}\times{\mathbb Z}_3$. The vertical axis stands for a
representation of ${\mathbb Z}_3$ as embedded in $[0,1)$. The
horizontal axis stands for the real line.  Since $d_{\beta}(1)=1030^\infty$,
there are two  complete central subtiles (green and red). The right figure
depicts  a zoom  along the vertical axis. This  zoom seems to suggest that the
central tile contains a full  stripe of the form
$[-\varepsilon,\varepsilon]\times {\mathbb Z}_3$, so that $\gamma(\beta)>0$. }
\label{fig:10_3} 
\end{figure}

\begin{figure}
\includegraphics[height=3cm, width=6cm]{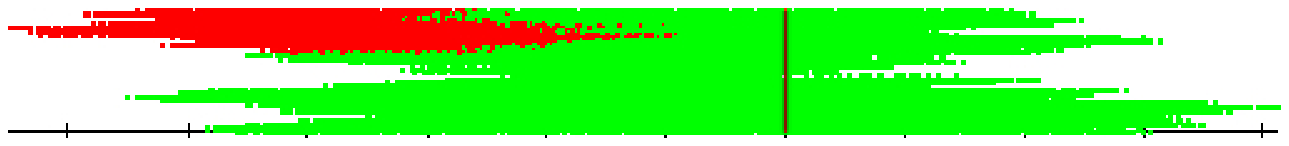}
\includegraphics[height=2.5cm]{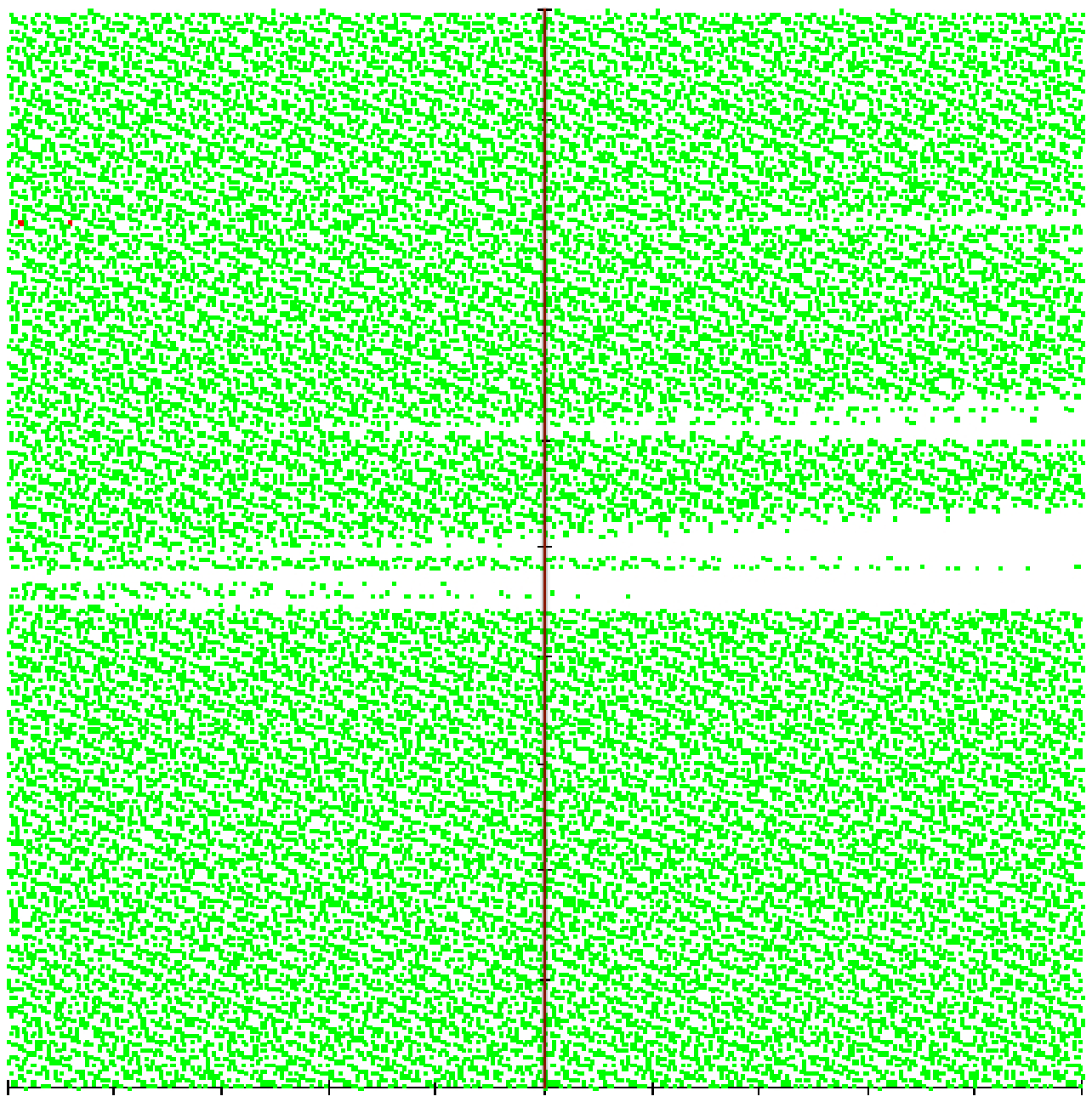}
\caption{A representation of the complete central tile for the Pisot number
  satisfying $\beta^2=4\beta + 3$. As in the previous case, we have
  $\GQ(\beta)=\GQ(\sqrt 7)$ and $N((\beta))=-3$. Thus the complete  central
  tile is again a subset of ${\mathbb R}\times{\mathbb Z}_3$. Since
  $d_{\beta}(1)=430^\infty$, there are two complete central subtiles (green and
red). The zoom  suggests that the  complete central tile  contains   
no stripe of the form  $[-\varepsilon,\varepsilon]\times {\mathbb Z}_3$,
so that $\gamma(\beta)=0$. }\label{fig:4_3}
\end{figure}

Let us recall (\cite{FrougnySolomyak92}, Proposition~1 and Lemma~3) that the
finitenes property  (F) holds  for any quadratic Pisot 
number $\beta$, and that those numbers are exactly the dominant root of the
polynomials $X^2-aX-b$ with  $a \geq b\geq 1$ or $q\geq 3$ and $-a+2\leq
b\leq -1$. Consequently, we may apply Theorem~\ref{thm:disjointness}, and
 the intersections between  complete  $x$-tiles determine  their boundary,
 which  have zero measure. The same property holds for the subtiles. We then
 use the fact that inner points of $x$-tiles and subtiles are exclusive to
 deduce an explicit formula for $\gamma(\beta)$.

\begin{theorem}\label{TheoGammaQuad}
If $\beta$ is quadratic, then $\gamma(\beta)$ is given by the
formula~\eqref{eq:minorationforgamma}, that is in that case an equality.
\end{theorem}

\begin{proof} First recall that since $\GK_\infty$ is one-dimensional, one has
  $\Delta(x)=x$ for all $x\in [0,1]$. We use the notation introduced in
  Proposition~\ref{prop:uvbounds}. We have to show that the lower bound is
  an upper bound too. We will show the following:\\
\begin{subequations}\label{eq:equality}
If $x \in\GZ[1/\beta]\cap (0,1)$ with $ -\Tt(x)\cap (0,1) \times 
  \overline{\phi_f(\GZ_{(N(\beta))})}\neq\emptyset$, then 
\begin{equation}\label{eq:equalityBax}
\gamma(\beta) \leq\inf\left\{ \pi_{\infty}\left(-\Tt(x) \cap (0,1) \times 
  \overline{\phi_f(\GZ_{(N(\beta))})}\right)\right\}.
\end{equation}
If $a\in\A$ with $-\Ta{a}\cap [T_\beta^{(a-1)},1)\times 
  \overline{\phi_f(\GZ_{(N(\beta))})}\neq\emptyset$, then
\begin{equation}\label{eq:equalityAab}
\gamma(\beta) \leq\inf\left\{ \pi_{\infty}\left(-\Ta{a} \cap (T_\beta^{(a-1)},1)
\times  \overline{\phi_f(\GZ_{(N(\beta))})}\right)\right\}.
\end{equation}
\end{subequations}
Since $-\Tt(x) \cap (0,1) \times 
  \overline{\phi_f(\GZ_{(N(\beta))})}\supset B_{a,x}$ for every $a\in\A$ and
  $-\Ta{a} \cap[T_\beta^{(a-1)},1) \times  \overline{\phi_f(\GZ_{(N(\beta))})}
  \supset A_{a,b}$ for every $b$ with $T_\beta^{(b-1)}(1)\leq
  T_\beta^{(a-1)}(1)$, the theorem will follow from~\eqref{eq:equality}. Notice
  that by continuity of $\pi_\infty$, taking open or closed intervals in
  $B_{a,x}$ or $A_{a,b}$ has no influence on the infimum we are interested in.\\

We begin with~\eqref{eq:equalityBax}.   Let $x\in\GZ[1/\beta] \cap
(0,1)$. Let  $z\in -\Tt(x)\cap (0,1)\times
\overline{\phi_f(\GZ_{(N(\beta))})}$. Since $\beta$ has degree 2, the property
(F) is satisfied, and $\Tt(x)$ is the closure of its subset of 
exclusive inner points by Proposition~\ref{prop:innerpoint}.  

Let us fix $\varepsilon>0$. There exists an exclusive inner point $y\in
-\Tt(x)\setminus ( -\Tt)$ such that $\Vert y-z\Vert\leq\varepsilon/2$.  
Since $y$ is an inner point and all inner points are exclusive, there exists
$\nu<\varepsilon/2$ such that the ball $B(y,\nu)$ is contained in $-\Tt(x)
\setminus (-\Tt)$. By Lemma~\ref{lem:localsurjective}, the set 
$\phb((\pi_{\infty}(y)-\nu, \pi_{\infty}(y)+\nu) \cap\GZ_{(N(\beta))})$ is
dense in $[\pi_{\infty}(y)-\nu,
\pi_{\infty}(y)+\nu]\times\overline{\phi_f(\GZ_{(N(\beta))})}$. Therefore, it
intersects $B(y,\nu)$, and there exists $w\in
(\pi_{\infty}(y)-\nu,\pi_{\infty}(y)+\nu) \cap\GZ_{(N(\beta))}$ 
such that $\phb(w)\in - \Tt(x)\setminus (-\Tt)$. For $w\leq
\pi_\infty(z)+\varepsilon$, we know by Theorem~\ref{thm:purpernonunit} that the
$\beta$-expansion of $w$ is not purely periodic. Hence $\gamma(\beta)\leq
\pi_{\infty}(z)+\varepsilon$. Finally, $\gamma(\beta)\leq\pi_{\infty}(z)$ and
\eqref{eq:equalityBax} is proved.\\

The proof for the upper bound~\eqref{eq:equalityAab} follows the same
lines. Let $z\in -\Ta{a}\cap (T_\beta^{(a-1)},1)\times
\overline{\phi_f(\GZ_{(N(\beta))})}$. Then $\Ta{a}$ is the closure of its  set
of exclusive inner points (with respect to $\Ta{b}$, $b\neq a$). For
$\varepsilon>0$, there exists an exclusive inner point $y$ and $\nu>0$ such
that $B(y,\nu)\subset-\Ta{a}\setminus\bigcup_{b\in\A\setminus\{ a\}}\Ta{b}$ and
$(\pi_{\infty}(y)-\nu,\pi_{\infty}(y)+\nu)\subset (T_\beta^{(a-1)},1)$ (this
second condition is the reason for which we take an open intervall
in~\eqref{eq:equalityAab}). By Lemma~\ref{lem:localsurjective}, there exists
$w\in (\pi_{\infty}(y)-\nu,\pi_{\infty}(y)+\nu) \cap\GZ_{(N(\beta))}$ such that
$\phi_\beta(w)\in -\Ta{a}\setminus\bigcup_{b\in\A\setminus\{ a\}}\Ta{b}$. Since
$\pi_\infty(w)>T_\beta^{(a-1)}$, $w\not\in\Pi_\beta^{(r)}$. Therefore,
$\gamma(\beta)\leq \pi_{\infty}(z)+\varepsilon$. Finally,
$\gamma(\beta)\leq\pi_{\infty}(z)$ and \eqref{eq:equalityAab} is proved.
\end{proof}

Suppose that  the degree of $\beta$ is  larger than 2.  We know that
$\pi_{\infty}(\phb({\mathbb Q} \cap [0,1])) \subset \Delta_{\infty}
([0,1])$. However, the  diagonal  set $\Delta_{\infty}([0,\infty))$  
 has empty interior in $\GK_\infty$. Consequently, it may happen that
 $\pi_{\infty}( -\Tt(x))$ is tangent to the diagonal $ \Delta_{\infty}
 ([0,\infty))$; in this latter case, $-\Tt(x)$  provides no point
 with a non-periodic beta-expansion and the conclusion of
Theorem~\ref{TheoGammaQuad} may fail.\\

\section{Two  quadratic examples}\label{sec:examples}

In the previous section, we have proved that $\gamma(\beta)$ is deeply related
with the intersections between subtiles and $x$-tiles.  In this section, we
will detail on two examples how $\gamma(\beta)$ can be explicitely computed. To
achieve this task, we will use the boundary graph defined in 
Section~\ref{sec:boundarygraph}. In Corollary~\ref{cor:boundarygraphalgorithm},
we have proved that the boundary graph can be computed by three conditions
(N1), (N3) and (N4). Conditions (N1) and (N4) are simple numerical 
conditions. On the contrary, condition (N3) implies the integer ring $\OO$. In
order to check this condition, we need to find an explicit basis of ${\mathfrak
  O} \cap {\mathbb Z}[1/\beta]$. We thus  introduce below  a sufficient
condition that reduces  ${\mathfrak  O} \cap {\mathbb Z}[1/\beta]$ to ${\mathbb
  Z}[\beta]$. 

\begin{lemma} \label{lem:quad}
Let $\beta$  be  such that $\beta{\mathfrak O}$ has only divisors of degree
1, and with inertia degree 1. Let $x \in\GZ[1/\beta]$. If $\beta^k x \in
{\mathfrak O}$, then $\beta^k x \in {\mathbb Z}[\beta]$. 
\end{lemma}

\begin{proof}
Let us expand $x$ as $x = a_{d-1} \beta^{d-1} + \dots + a_0 + \dots +
a_{-N}\beta^{-N}$, with $a_i \in {\mathbb Z}$ (it is not the
$\beta$-expansion).  If $N > k$, then $\beta^N x = \beta^{N-k} (\beta^k x) \in
\beta^{N-k}{\mathfrak O}$.  We deduce that $a_{-N} \in \beta^{N-k}{\mathfrak O}
+ \beta {\mathbb Z}[\beta] \subset \beta{\mathfrak O}$. Hence $a_{-N} \in
\beta{\mathfrak O} \cap\GZ$. Since $\beta{\mathfrak O}$ has only divisors of
degree 1 and with inertia degree 1,  $N(\beta)$ divides $a_{-N}$. From
$N(\beta)/\beta \in {\mathbb Z}[\beta]$, we deduce that $a_{-N}/\beta \in
{\mathbb Z}[\beta]$. Then $x$ admits an expansion of size at most
$\beta^{-N+1}$: $x = b_{d-1} \beta^{d-1}+\cdots  b_0+\cdots
+b_{-N+1}\beta^{-N+1}$. We conclude by induction that $\beta^k x\in\GZ[\beta]$.  
\end{proof}

Let us stress the fact that if $\beta$ is a quadratic number that satisfies the
conditions of Proposition ~\ref{prop:charactdensity}, then Lemma~\ref{lem:quad}
holds.  In this case, Corollary~\ref{cor:boundarygraphalgorithm} reads as
follows to compute the boundary graph. 

\begin{corollary}
Suppose that $\beta$ is a quadratic number such that  $\beta{\mathfrak O}$ has
only divisors of degree 1 and inertia degree 1. Let $\beta^2 = a \beta + b$ be 
its minimal polynomial. The boundary graph of $\beta$ can be explicitely
computed as follows. 
\begin{enumerate}
\item Consider all triplets $[a,x,b]$ such that $x=K+\beta L$, $(K,L) \in
  {\mathbb Z}^2$, with 
\begin{itemize}
\item  $K \leq \frac{\beta-a+3a\beta-\beta^2-a^2}{(2\beta-a)(1+a-\beta)}$ and
  $L \leq \frac{1+2a - \beta}{(2\beta-a)(1+a-\beta)}$.  
\item  $-T_{\beta}^{(a-1)}(1)<x< T_{\beta}^{(b-1)}(1)$ and $a\not=b$ if $x=0$.
\end{itemize}
\item Put an edge betwenn two triplets $[a,x,b]$ and $[a_1,x_1,b_1]$ if there
  exists $q_1$ and $p_1$ such that 
\begin{itemize}
\item $x_1 = \beta^{-1}(x+q_1 - p_1)$,
\item $a_1 \xrightarrow{p_1}a$ and $b_1 \xrightarrow{q_1} b$ are edges of
  the admissibility graph. 
\end{itemize}
\item Recursively remove edges that have no outgoing edge. 
\end{enumerate}
\end{corollary}

\begin{proof}
From the proof of Corollary \ref{cor:boundarygraphalgorithm},
it is sufficient to exhibit a set that contains  all the triplets $[a,x,b]$
satisfying conditions (N1), (N3) and (N4). Then the recursive deletion of edges
will reduce the graph to the exact boundary 
graph. In this case, condition (N3) implies that $x \in {\mathbb Z}[\beta]$. 
Then we are looking for all $x$'s such that $x= K \beta + L$, with $K,L \in
{\mathbb Z}$, and such that conditions (N1) and (N4) are satisfied.   
let $x_2$ denote the conjugate of $x$ and $\beta_2 = a-\beta $ denote
the conjugate of $\beta$. We obtain $K = (\beta_2 x-\beta x_2)/(\beta-\beta_2)$
and $L =(x-x_2)/(\beta-\beta_2)$.  
Condition (N1) means that $x \leq 1$, and condition (N4) implies that $|x_2|
\leq \frac{\lfloor\beta\rfloor}{1-|\beta_2|} = \frac{a}{1+a-\beta}$.  
We deduce that if  $[a,x,b]$  satisfies 
the three conditions (N1), (N3) and (N4), then $x = K \beta + L$ with $K \leq
\frac{\beta-a+3a\beta-\beta^2-a^2}{(2\beta-a)(1+a-\beta)}$  
and $L \leq \frac{1+2a - \beta}{(2\beta-a)(1+a-\beta)}$. 
\end{proof}

\begin{figure}[t]
\includegraphics[width=12cm]{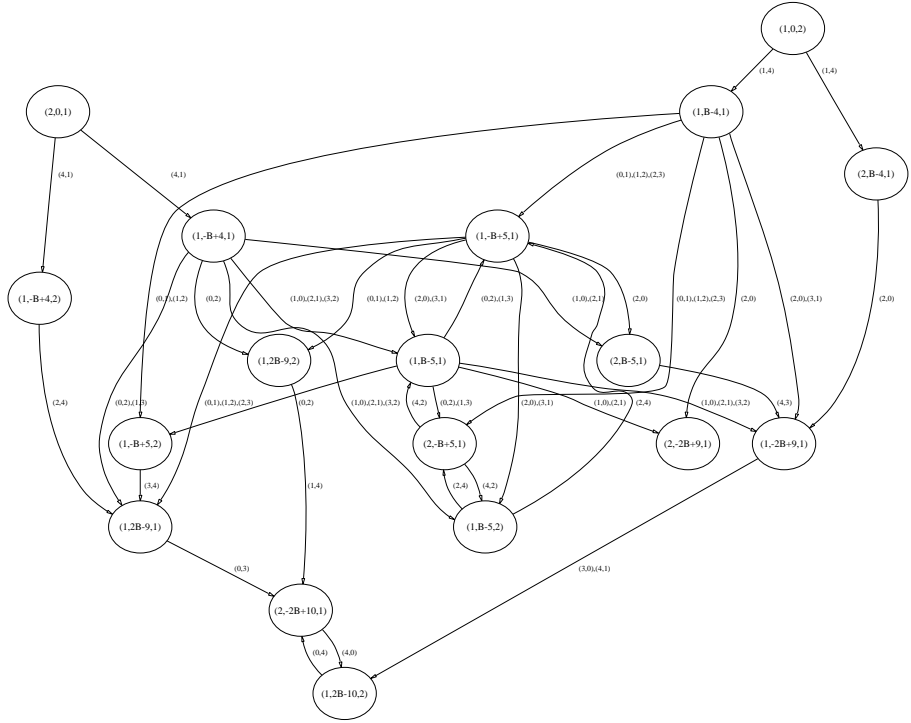}
\caption{ The boundary graph for $\beta^2=4\beta+3$.}\label{BG1} 
\end{figure}

\begin{figure}[t]
\begin{center}
\includegraphics[width=8cm]{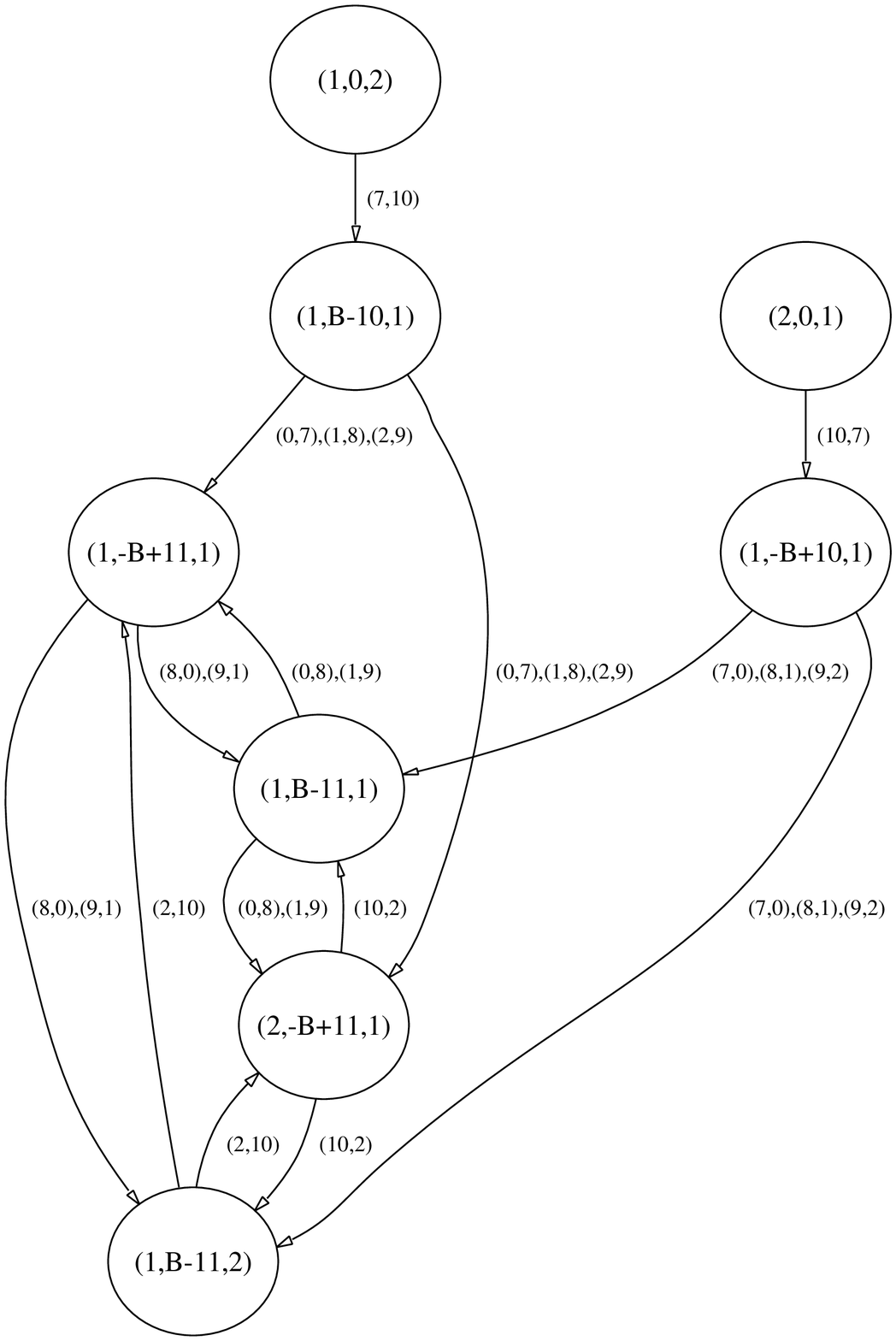}
\includegraphics[width=7cm]{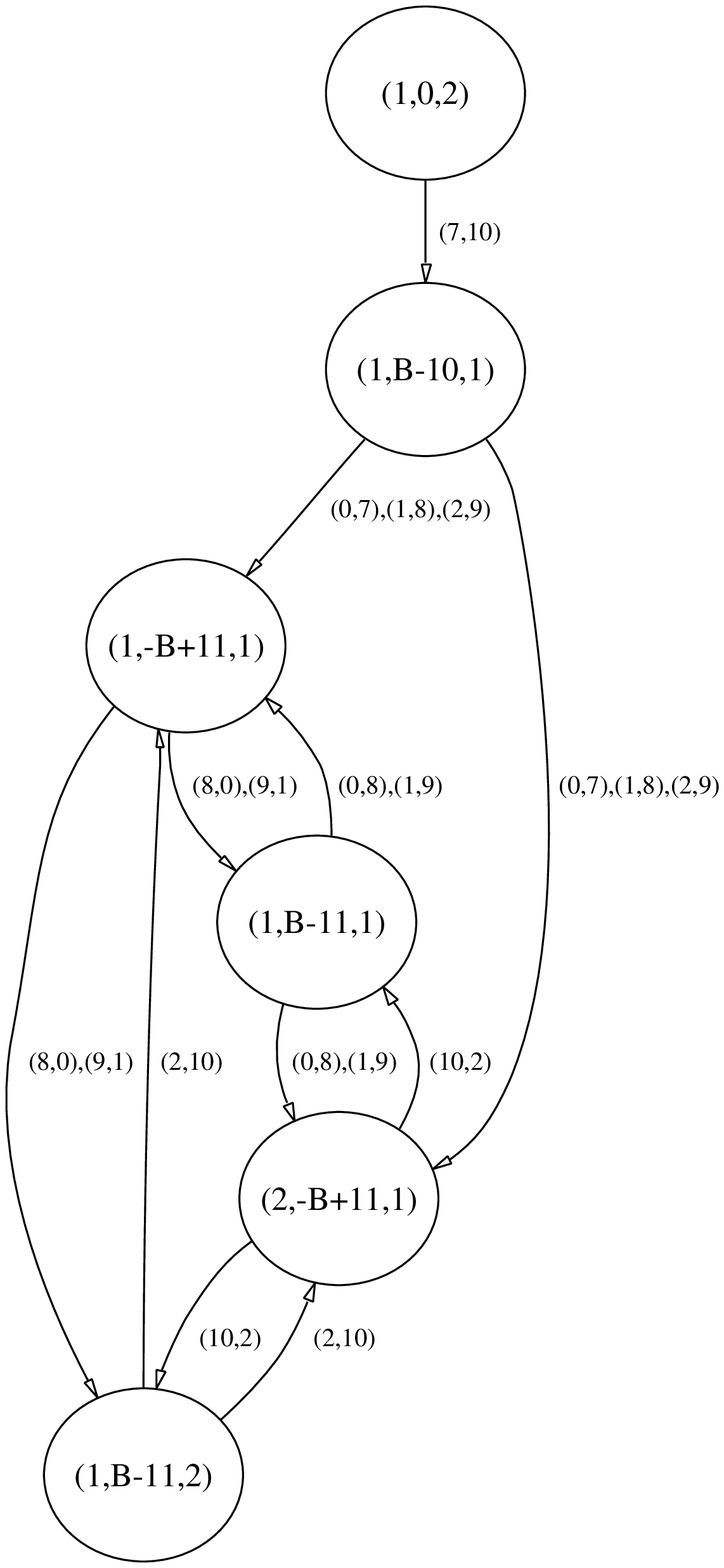}
\end{center}
\caption{ (Left) The boundary graph for $\beta^2=10\beta+3$.  The notation $B$ stands for $\beta$.   \newline (Right)
Subgraph of the boundary graph  that gathers
infinite paths starting from a node $[a,x,b]$ with $x \in {\mathbb Z}[1/\beta]$, $0 \leq x < T^{(b-1)}(1)$  and $a<b$ if $x=0$. These nodes are $[1,0,2]$,  $[1,\beta-10,1]$, $[1,11-\beta,1]$ and $[2,11-\beta,1]$.  They exacly stand for the set of intersections that contribute to the computation of $\gamma(\beta)$ in Proposition \ref{prop:uvbounds}.  
}\label{BG2Bis}
\end{figure}

When $\beta^2=4\beta+3$ the bounds are $K \leq 11$ and $L \leq 3$.  We deduce
that the boundary graph contains 18 nodes (Fig. \ref{BG1}). 
If $[a,x,b]$ is a node of the boundary graph, we have $x \in \pm \{0,\beta-4,
5-\beta,2\beta-10,2\beta-9\}$. When  $\beta>1$  defined by
$\beta^2=10\beta+3$,  the bounds are $K \leq 14$ and $L \leq 2$. The boundary
graph contains 8 nodes (see left side of Fig. \ref{BG2Bis}). If $[a,x,b]$ is a
node of the boundary graph, we have $x \in\pm \{0,11-\beta,\beta-10\}$.  

\begin{figure}[t]
\includegraphics[width=14cm]{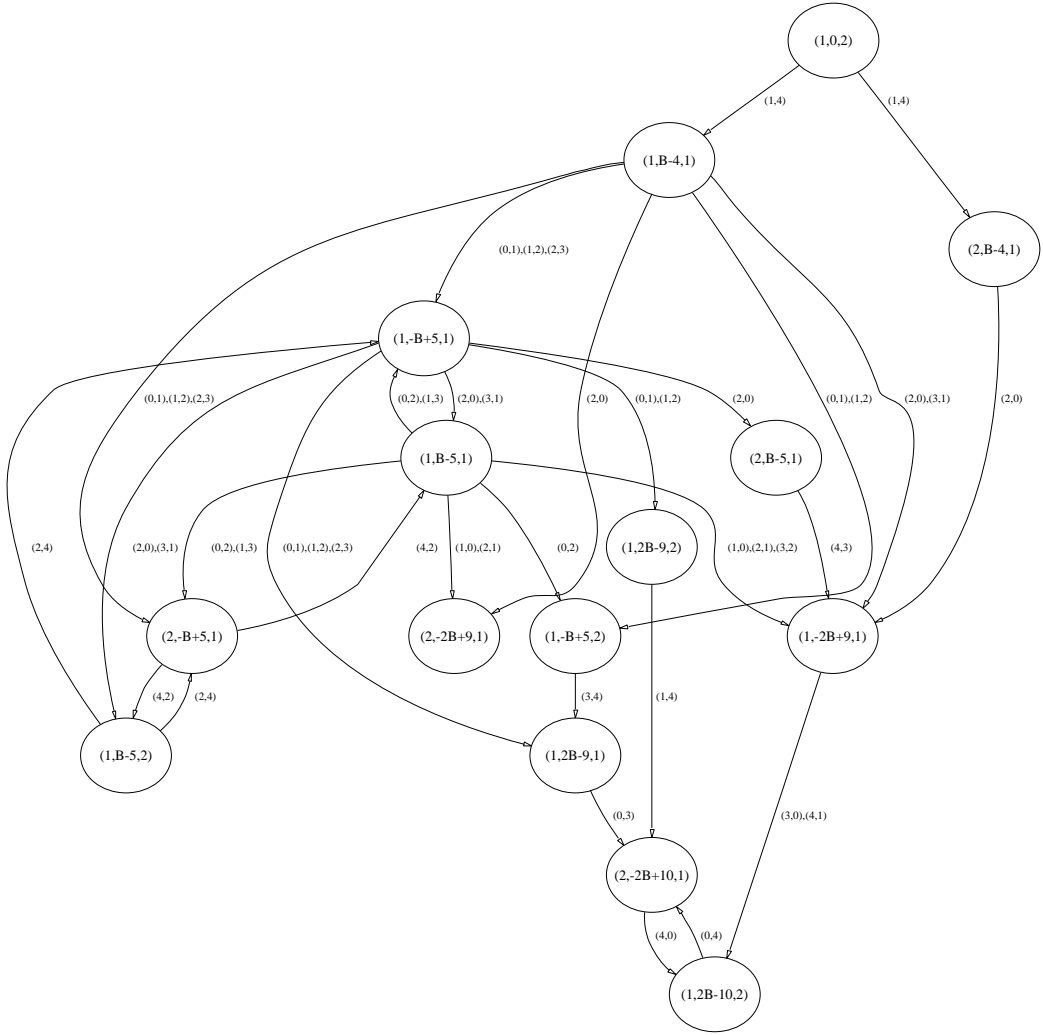}
\caption{   Subgraph of the boundary graph for $\beta^2=4\beta+3$ that gather
infinite paths starting from a node $[a,x,b]$ with $x \in {\mathbb Z}[1/\beta]$, $0 \leq x < T^{(b-1)}(1)$  and $a<b$ if $x=0$. These nodes are  $[1,0,2]$, $[1,2\beta-9,1]$,
$[1,2\beta-9,2]$, $[1,\beta-4,1]$, $[2,\beta-4,1]$, $[1,-\beta+5,1]$,
$[1,-\beta+5,2]$, $[2,-\beta+5,1]$, $[2, -2\beta+10,1]$.  They exacly stand for the set of intersections  that  contribute to the computation of $\gamma(\beta)$ in Proposition \ref{prop:uvbounds}.  }\label{BG1Bis} 
\end{figure}

\begin{proposition}
Let $\beta>1$  defined by $\beta^2=4\beta+3$. 
 There are 9 non-empty intersections between the  central subtiles and
the  neighbouring $x$-tiles, namely  $\widetilde{\mathcal T}^{(1)} \cap
\widetilde{\mathcal T}^{(2)}$,  
 $\widetilde{\mathcal T}^{(1)} \cap (\widetilde{\mathcal T}^{(1)} + \phb(2\beta
 -9))$, 
$\widetilde{\mathcal T}^{(1)} \cap (\widetilde{\mathcal T}^{(2)} + \phb(2\beta
-9))$, $\widetilde{\mathcal T}^{(1)} \cap (\widetilde{\mathcal T}^{(1)} +
\phb(\beta-4))$, $\widetilde{\mathcal T}^{(2)} \cap (\widetilde{\mathcal
  T}^{(1)} + \phb(\beta-4))$, $\widetilde{\mathcal T}^{(1)} \cap
(\widetilde{\mathcal T}^{(1)} + \phb(5-\beta))$,  $\widetilde{\mathcal T}^{(2)}
\cap (\widetilde{\mathcal T}^{(1)} + \phb(5-\beta))$, $\widetilde{\mathcal T}^{(1)}
\cap (\widetilde{\mathcal T}^{(2)} + \phb(5-\beta))$  $\widetilde{\mathcal
  T}^{(2)} \cap (\widetilde{\mathcal T}^{(1)} + \phb(10-2\beta))$. The
expansions of the points lying in one of those intersections are constrained by
the graph depicted in  Fig.~\ref{BG1Bis}. 
\end{proposition}

\begin{proof}
In order to obtain the interesting intersections, we consider in the boundary
graph the subgraph of paths starting from  
$[a,x,b]$ with $x \in [0, T_{\beta}^{(b-1)}(1)[$ and $a<b$ if $x=0$. In the
boundary graph, there are  9  nodes which satisfy 
these conditions: $[1,0,2]$, $[1,2\beta-9,1]$,
$[1,2\beta-9,2]$, $[1,\beta-4,1]$, $[2,\beta-4,1]$, $[1,-\beta+5,1]$,
$[1,-\beta+5,2]$, $[2,-\beta+5,1]$, $[2, -2\beta+10,1]$. From these nodes, infinite paths 
span a subgraph with 15 nodes, depicted in  Fig. \ref{BG1Bis}.
\end{proof}


We obtain another graph for $\beta^2=10\beta+3$.

\begin{proposition}
Let $\beta>1$  defined by  $\beta^2=10\beta+3$. There are exactly 4 non-empty
intersection between the  central subtiles and $x$-tiles, namely  
$\widetilde{\mathcal T}^{(1)} \cap \widetilde{\mathcal T}^{(2)}$,
$\widetilde{\mathcal T}^{(1)} \cap (\widetilde{\mathcal T}^{(1)} + \phb(\beta
-10))$, 
 $\widetilde{\mathcal T}^{(1)} \cap (\widetilde{\mathcal T}^{(1)} + \phb(-\beta
 +11))$, $\widetilde{\mathcal T}^{(2)} \cap (\widetilde{\mathcal T}^{(1)} +
 \phb(-\beta +11))$. The expansions of the points lying in one of those
 intersections are constrained by the graph depicted in
 Fig. \ref{BG2Bis}. 
\end{proposition} 
\begin{proof}
In the boundary graph, nodes $[a,x,b]$ that satisfy the condition  $x \in [0,
T_{\beta}^{(b-1)}(1)[$ and $a<b$ if $x=0$ are 
$[1,0,2]$,  $[1,\beta-10,1]$, $[1,11-\beta,1]$ and $[2,11-\beta,1]$. In the
boundary graph, paths starting from these nodes cover 
a subgraph with 5 nodes, shown in    Fig. \ref{BG2Bis}.
\end{proof}

We now have the tools to compute $\gamma(\beta)$ in some specific cases.

\begin{lemma}\label{lemmaBG1}
Let $\beta^2 = 4\beta+3$. We recall that $\pi_{\infty}$ stands for the
projection from ${\mathbb K}_{\infty}$ to ${\mathbb R}$. Then 
$$\pi_{\infty}\left(\widetilde{\mathcal T}^{(1)}\cap( \widetilde{\mathcal
    T}^{(1)}+\phb(\beta-4)) \right)  =\left\{ \sum a_{i} {\beta_2}^{i} 
\ |\ 
a_{2i}\in \{0,1,2\}, a_{2i+1} \in \{2,3,4\}
\right\}.$$
\end{lemma}
\begin{proof}
We use the boundary subgraph depicted in Fig.~\ref{BG1Bis}. 
By  construction, any point of the intersection ${\widetilde{\mathcal
    T}}^{(1)}\cap\left( \widetilde{\mathcal T}^{(1)}+\phb(\beta-4) \right)$ 
can be expanded as $z = \sum_{i \geq 0} p_i \phb(\beta^i)$, where $(p_i)$ is
the first coordinate of the labeling of a path starting  in $[1,\beta-4,1]$. By
looking at paths starting from $[1,\beta-4,1]$ we  check in the graph that such
sequences $(p_i)$'s satisfy   $p_{2i} \in \{0,1,2\}$ and $p_{2i+1} \in
\{2,3,4\}$. Conversely, we also check that every sequence of this form is the
first coordinate of the labeling of a path starting in $[1,\beta-4,1]$ in the
graph. This yields 
$$\pi_{\infty}(\widetilde{\mathcal T}^{(1)}\cap(\widetilde{\mathcal
  T}^{(1)}+\phb(\beta-4)))  = \left\{ \sum_{i\geq 0} \{0,1,2\} \beta_2^{2i} +
  \sum_{i\geq 0} \{2,3,4\} \beta_2^{2i+1} \right\} .$$ 
\end{proof}

In order to compute $\gamma(\beta)$, we use  the following folklore lemma.
\begin{lemma}[Cookie Cantor Lemma] Let $\alpha<1$ be an integer number 
$$
X(\alpha,n):= \left\{\left. \sum_{i\geq 0} a_i \alpha^i \right| a_i \in
  \{0,1,\dots,n-1\} \right\} \subset \left[0, \frac {n-1}{1-\alpha}\right]
$$
The two end points $\{0, \frac{n-1}{1-\alpha}\}$  belong to $X(\alpha,n)$.
Furthermore,  if $\alpha>1/n$, then it is a Cantor cookie cutter set and if 
$\alpha\in [1/n,1)$, then $X(\alpha,n)$ coincides with the interval 
$\left[0, (n-1)/(1-\alpha)\right]$.
\end{lemma}

\begin{proof}
The set $X(\alpha,n)$ is the attractor of the IFS:
$X= \bigcup_{i=0}^{n-1} \alpha X+i$
which has a unique non-empty compact solution. 
It is easy to see that the right hand side
is a solution if $\alpha\in [1/n,1)$.
\end{proof}

\begin{theorem}\label{thm:2+sqrt7}
One has 
$$ \gamma(2+\sqrt 7) = 0.$$
\end{theorem}
\begin{proof}
If $\beta^2 = 4\beta+3$ and $\beta$ is  a Pisot number, then $\beta = 2+\sqrt
7$.  We also check that $\beta$ satisfies the conditions of Corollary \ref{cor:quadraticdensity},
hence $\overline{\phi_f(\GZ_{(N(\beta))})}={\mathbb Z}_3$. In this case, the set $A_{a,b}$ and
$B_{a,x}$ in Proposition \ref{prop:uvbounds} simply correspond to intersections between tiles,
with no more diagonal set: $A_{a,b}=\pi_\infty(-\Ta{a}\cap -\Ta{b} \cap [T_\beta^{(b-1)}(1), T_\beta^{(a-1)}(1)])\times {\mathbb Z}_3)$
and $B_{a,x}=\pi_\infty(-\Ta{a}\cap -\Tt(x)\cap (0, T_\beta^{(a-1)}(1))\times {\mathbb Z}_3)$. 
Then computing $\gamma(\beta)$ reduced to understanding intersections between tiles.

Let $-\alpha$ denote the conjugate of $\beta$, that is,
$\alpha=\sqrt{7}-2>1/3$. Lemma~\ref{lemmaBG1} exhibits a set that we need to
compute explicitly.  
\begin{eqnarray*}
\pi_{\infty}\left(\widetilde{\mathcal T}^{(1)}\cap( \widetilde{\mathcal
    T}^{(1)}+\phb(\beta-4)) \right)&=& \left\{ 
\sum a_{i} {\alpha}^{2i} - \sum b_{i} {\alpha}^{2i+1} 
\ |\ 
a_i\in \{0,1,2\}, b_i \in \{2,3,4\}
\right\} \\
&=&
\left\{
\sum a_{i} \alpha^{2i} +\sum c_{i} \alpha^{2i+1} 
\ |\ 
a_i\in \{0,1,2\} , c_i\in  \{-2,-3,-4\}
\right\}\\
&=& X(\alpha,3) - \alpha (4+ 4\alpha^2+4\alpha^4+\dots )\\
&=& \left[ -\frac{4\alpha}{1-\alpha^2} ,  -\frac{4\alpha}{1-\alpha^2}+ \frac
  2{1-\alpha}  \right]\\ 
&=& \left[  -\frac{4\alpha}{1-\alpha^2} ,  \frac{2 - 2 \alpha}{1-\alpha^2}
\right] \ni 0 
\end{eqnarray*}
Hence zero is the minimum of
$[0,T_{\beta}(1)]\cap\pi_{\infty}(-\widetilde{\mathcal
  T}^{(1)}\cap(-\widetilde{\mathcal T}(\beta-4))$ 
and Theorem \ref{TheoGammaQuad} implies that $\gamma(\beta)=0$. 
\end{proof}

A completely different behaviour appears when modifying only one digit
in the quadratic equation satisfied by $\beta$.

\begin{theorem}\label{thm:5+2sqrt7}
One has
$$\gamma(5+2\sqrt 7) = \frac{7 - \sqrt7}{12} $$
\end{theorem}

\begin{proof}

The number $5+2\sqrt 7$ is the root of $\beta^2-10\beta-3=0$. As before,  conditions of Corollary \ref{cor:quadraticdensity} are satisfied
hence $\overline{\phi_f(\GZ_{(N(\beta))})}={\mathbb Z}_3$, and studying  intersections between tiles is enough to
compute $\gamma(\beta)$.

From the graph depicted  in Fig. \ref{BG2Bis}, we deduce that non-empty
intersections in the numeration tiling are given by 
$ \widetilde{\mathcal T}^{(1)}\cap(\widetilde{\mathcal T}^{(2)}$, 
$ \widetilde{\mathcal T}^{(1)}\cap(\widetilde{\mathcal
  T}^{(1)}+\phb(11-\beta))$,  
$ \widetilde{\mathcal T}^{(2)}\cap(\widetilde{\mathcal
  T}^{(1)}+\phb(11-\beta))$,  
and  $ \widetilde{\mathcal T}^{(1)}\cap(\widetilde{\mathcal
  T}^{(1)}+\phb(+10+\beta))$. 

We can detail the expansion of the real projection of the three last sets.

$$ \pi_{\infty}(\Ta{1}\cap(\Ta{1}+\phb(11-\beta))) = \{8,9\} + \beta_2
\sum_{i\geq 0} \{0,1,2\}\beta_2^{2i} +  \{8,9,10\}\beta_2^{2i+1} $$  

$$ \pi_{\infty}(\widetilde{\mathcal T}^{(2)}\cap(\widetilde{\mathcal T}^{(1)}+\phb(11-\beta))) = 10 + \beta_2  \sum_{i\geq 0} \{0;1,2\}\beta_2^{2i} +  \{8,9,10\}\beta_2^{2i+1} $$ 

$$ \pi_{\infty}(\widetilde{\mathcal T}^{(1)}\cap(\widetilde{\mathcal
  T}^{(1)}+\phb(-10+\beta))) =  \{0,1,2\} + \beta_2 \sum_{i\geq 0}
\{8,9,10\}\beta_2^{2i} +  \{0,1,2\}\beta_2^{2i+1} $$

We use the Cookie Cantor Lemma stated above with $\alpha = -\beta_{(2)}=
3\beta^{-1}$  and $n=3> \alpha^{-1}$ to compute the sum that is involved in
each intersection. 
\begin{eqnarray*}
\sum_{i\geq 1} \{0,1,2\}\beta_2^{2i} +  \{8,9,10\}\beta_2^{2i+1}  &= &
\sum_{i\geq 0} \{0;1,2\}\alpha^{2i} -  \{8,9,10\}\alpha^{2i+1} \\ 
&= & -10 \sum_{i \geq 0}\alpha^{2i+1}  + \sum_{j\geq 0} \{0,1,2\} \alpha^j \\
& = & \frac{-10 \alpha}{1 - \alpha^2} + \left[ 0, \frac{2}{1-\alpha} \right]  =
\left[\frac{-10 \alpha}{1 - \alpha^2}, \frac{-8 \alpha + 2}{1 - \alpha^2}
\right]  
\end{eqnarray*}

Similarly, we have
\begin{eqnarray*}
\sum_{i\geq 1} \{8,9,10\}\beta_2^{2i} +  \{0,1,2\}\beta_2^{2i+1}  &= & 10
\sum_{i \geq 0}\alpha^{2i}  - \sum_{j\geq 0} \{0,1,2\} \alpha^j \\ 
& = & \frac{10 }{1 - \alpha^2} - \left[ 0, \frac{2}{1-\alpha} \right]  =
\left[\frac{8 - 2 \alpha}{1 - \alpha^2}, \frac{10}{1 - \alpha^2} \right]. 
\end{eqnarray*}

We deduce that
 $$ \pi_{\infty}(\widetilde{\mathcal T}^{(1)}\cap(\widetilde{\mathcal
   T}^{(1)}+\phb(11-\beta))) = \left[ 8 -\alpha \frac{-8 \alpha + 2}{1 -
     \alpha^2}, 8 - \alpha \frac{- 10 \alpha }{1 - \alpha^2} \right] \subset
 \left] 0, \infty \right[ .$$ 

Hence  $-\pi_{\infty}(\widetilde{\mathcal T}\cap(\widetilde{\mathcal
  T}^{(1)}+\phb(11 - \beta))) \cap [0,1] = \emptyset$. Similarly, 
we have  $-\pi_{\infty}(\widetilde{\mathcal T}\cap(\widetilde{{\mathcal
    T}^{(2)}}+\phb(11 - \beta))) \cap [0,1] = \emptyset$, 
 so that both  intersections cannot be taken into account in the computation of
 $\gamma(\beta)$. This implies  that $-\pi_{\infty}(\widetilde{\mathcal
   T}^{(2)}) \cap [0,\infty]$ does not intersect the projection of any tile
 $-\pi_{\infty}(\widetilde{\mathcal T}(x))$.

We also have
$$\pi_{\infty}(\widetilde{\mathcal T}^{(1)}\cap(\widetilde{\mathcal
  T}^{(1)}+\phb(\beta-10))) =  \left[\frac{-10 \alpha}{1 - \alpha^2}, \frac{-8
    \alpha + 2}{1 - \alpha^2} \right].$$  

Hence, the minimum of $-\pi_{\infty}(\widetilde{\mathcal
  T}^{(1)}\cap(\widetilde{\mathcal T}^{(1)}+\phb(\beta-10)))$ 
is $\frac{8 \alpha - 2}{1 - \alpha^2}$.

In order to apply Theorem \ref{TheoGammaQuad},  we  prove that the infimum of intersections of the form $\pi_{\infty}(A_{a,b})$ (situation (A) or (B)) is strictly larger than the infimum of  intersections  $\pi_{\infty}(B_{x,a})$ (situation (C)). 
By  definition, we have   
$\pi_{\infty}(\widetilde{\mathcal T}^{(2)})= \{ \sum_{i \geq 0} a_{2i}
\alpha^2i - \sum_{i \geq 0} a_{2i+1} \alpha^{2i+1}\}$ 
where sequences $a_1 \dots a_i \dots$ are sequences starting from 2 in the
reverse of the admissibility graph. We deduce 
 that $a_0 = 10$, $a_1 \leq 9$, $a_2 \geq 0$, and then, $a_{2i+2} \geq 0$ and
 $a_{2i+3} \leq 10$. Hence 
$$\min \pi_{\infty}(\widetilde{\mathcal T}^{(2)} )\geq 10 - 9 \alpha + 0 \alpha^2 - 10
\alpha^3 + \dots = 10 - 9 \alpha + 10 \frac{\alpha^3}{1-\alpha^2} >0.$$ 

Consequently,  $-\pi_{\infty}(\widetilde{\mathcal T}^{(2)}) \cap [0,\infty] =
\emptyset$ and situations (A) or (B) do not contribute to $\gamma(\beta)$.

From Theorem \ref{TheoGammaQuad}, we deduce that $\gamma(\beta)= \min  -
\pi_{\infty}(\widetilde{\mathcal T}^{(1)}\cap(\widetilde{\mathcal
  T}^{(1)}+\phb(11-\beta))) =  \frac{8\alpha - 2}{1-\alpha^2}= \frac{7 -
  \sqrt7}{12}$.  
\end{proof}

\section{Perspectives}

At least two  main directions deserve now to be discussed.  In the quadratic
case, what is the structure of the intersection graph allowing to compute
$\gamma(\beta)$? The first question is 
whether we can obtain an algorithmic way to compute $\gamma(\beta)$  for every
quadratic $\beta$.  
Then, can we deduce a general formula for $\gamma(\beta)$ for subfamilies of
$\beta$? The first step would be to describe properly the structure of the
boundary graph, at least for the family $\beta^2 = n \beta + 3$.\\ 

Another direction lies in the application of these methods in the three (or
more dimensional case), including 
the unit case.   At the moment we cannot give an explicit formula for
$\gamma(\beta)$. In order to generalise the results 
to higher degrees, an approximation of exclusive inner points by the diagonal
line of ${\mathbb K}_{\beta}$ is needed. This 
seems reasonable at least in the unit case, but requires a precise study of the
topology of the central tile.  Examples of computations 
of intersections between line and fractals are obtained in \cite{Akiyama05},
by numeric approximations.  As an example, an intersection is prove to be approximated by 0.66666666608644067488. 
Then it is not equal to 2/3, though very near from it. 
Theorem~\ref{thm:5+2sqrt7} is an example where we were able to compute
explicitely the value of $\gamma(\beta)$ and it turned out that
$\gamma(\beta)\in\GQ(\beta)$. The question of the algebraic nature of
$\gamma(\beta)$ in general is interesting.

\newpage

\footnotesize{
\bibliography{ABBS_Submitted}
\bibliographystyle{alpha}
}

\end{document}

%% file: graph2.pstex_t
\begin{picture}(0,0)%
\includegraphics{graph2.pstex}%
\end{picture}%
\setlength{\unitlength}{3947sp}%
\begingroup\makeatletter\ifx\SetFigFont\undefined%
\gdef\SetFigFont#1#2#3#4#5{%
  \reset@font\fontsize{#1}{#2pt}%
  \fontfamily{#3}\fontseries{#4}\fontshape{#5}%
  \selectfont}%
\fi\endgroup%
\begin{picture}(3501,2150)(8561,-4148)
\put(9026,-2259){\makebox(0,0)[b]{\smash{{\SetFigFont{5}{6.0}{\familydefault}{\mddefault}{\updefault}$t_1$}}}}
\put(11641,-2281){\makebox(0,0)[b]{\smash{{\SetFigFont{5}{6.0}{\familydefault}{\mddefault}{\updefault}$t_{n-1}$}}}}
\put(10954,-2280){\makebox(0,0)[b]{\smash{{\SetFigFont{5}{6.0}{\familydefault}{\mddefault}{\updefault}$t_{m+1}$}}}}
\put(11210,-2058){\makebox(0,0)[lb]{\smash{{\SetFigFont{5}{6.0}{\familydefault}{\mddefault}{\updefault}$t_{n}$}}}}
\put(9916,-2250){\makebox(0,0)[b]{\smash{{\SetFigFont{5}{6.0}{\familydefault}{\mddefault}{\updefault}$t_{m-1}$}}}}
\put(10386,-2250){\makebox(0,0)[b]{\smash{{\SetFigFont{5}{6.0}{\familydefault}{\mddefault}{\updefault}$t_{m}$}}}}
\put(9457,-2938){\makebox(0,0)[b]{\smash{{\SetFigFont{5}{6.0}{\familydefault}{\mddefault}{\updefault}$0,\cdots,t_m-1$}}}}
\put(9884,-3642){\makebox(0,0)[b]{\smash{{\SetFigFont{5}{6.0}{\familydefault}{\mddefault}{\updefault}$0,\cdots,t_{m+2}-1$}}}}
\put(10360,-4004){\makebox(0,0)[b]{\smash{{\SetFigFont{5}{6.0}{\familydefault}{\mddefault}{\updefault}$0,\cdots,t_{n}-1$}}}}
\put(8991,-2691){\makebox(0,0)[b]{\smash{{\SetFigFont{5}{6.0}{\familydefault}{\mddefault}{\updefault}$0,\cdots,$}}}}
\put(9005,-2790){\makebox(0,0)[b]{\smash{{\SetFigFont{5}{6.0}{\familydefault}{\mddefault}{\updefault}$t_2-1$}}}}
\put(9685,-3286){\makebox(0,0)[b]{\smash{{\SetFigFont{5}{6.0}{\familydefault}{\mddefault}{\updefault}$0,\cdots,t_{m+1}-1$}}}}
\put(8700,-2446){\makebox(0,0)[rb]{\smash{{\SetFigFont{5}{6.0}{\familydefault}{\mddefault}{\updefault}$0,\cdots,t_1-1$}}}}
\put(11045,-2571){\makebox(0,0)[b]{\smash{{\SetFigFont{8}{9.6}{\familydefault}{\mddefault}{\updefault}$m+2$}}}}
\put(11970,-2566){\makebox(0,0)[b]{\smash{{\SetFigFont{8}{9.6}{\familydefault}{\mddefault}{\updefault}$n$}}}}
\put(10598,-2571){\makebox(0,0)[b]{\smash{{\SetFigFont{8}{9.6}{\familydefault}{\mddefault}{\updefault}$m+1$}}}}
\put(10155,-2565){\makebox(0,0)[b]{\smash{{\SetFigFont{8}{9.6}{\familydefault}{\mddefault}{\updefault}$m$}}}}
\put(8784,-2575){\makebox(0,0)[b]{\smash{{\SetFigFont{8}{9.6}{\familydefault}{\mddefault}{\updefault}1}}}}
\put(9251,-2575){\makebox(0,0)[b]{\smash{{\SetFigFont{8}{9.6}{\familydefault}{\mddefault}{\updefault}2}}}}
\end{picture}%

%% file: ABBS_Submitted.bbl
\begin{thebibliography}{BFGK98}

\bibitem[Aki98]{Akiyama98}
S.~Akiyama.
\newblock Pisot numbers and greedy algorithm.
\newblock In {\em Number theory (Eger, 1996)}, pages 9--21. de Gruyter, Berlin,
  1998.

\bibitem[Aki00]{Akiyama00}
S.~Akiyama.
\newblock Cubic {P}isot units with finite beta expansions.
\newblock In {\em Algebraic number theory and Diophantine analysis (Graz,
  1998)}, pages 11--26. de Gruyter, Berlin, 2000.

\bibitem[Aki02]{Akiyama02}
S.~Akiyama.
\newblock On the boundary of self affine tilings generated by {P}isot numbers.
\newblock {\em J. Math. Soc. Japan}, 54(2):283--308, 2002.

\bibitem[Aki07]{Akiyama07}
S.~Akiyama.
\newblock Pisot number system and its dual tiling.
\newblock In {\em Physics and Theoretical Computer Science (Cargese, 2006)},
  pages 133--154. IOS Press, 2007.

\bibitem[AS05]{Akiyama05}
S.~Akiyama and K.~Scheicher.
\newblock Intersecting two dimensional fractals and lines.
\newblock {\em Acta Sci. Math. (Szeged)}, 3-4:555--580, 2005.

\bibitem[BBLT06]{Baratetal}
G.~Barat, V.~Berth{\'e}, P.~Liardet, and J.~M. Thuswaldner.
\newblock Dynamical directions in numeration.
\newblock {\em Ann. Inst. Fourier (Grenoble)}, 56(7):1987--2092, 2006.

\bibitem[Ber77]{BM77}
A.~Bertrand.
\newblock D\'eveloppements en base de {P}isot et r\'epartition modulo {$1$}.
\newblock {\em C. R. Acad. Sci. Paris S\'er. A-B}, 285(6):A419--A421, 1977.

\bibitem[BFGK98]{BFGK}
{\v{C}}.~Burd{\'{\i}}k, Ch. Frougny, J.~P. Gazeau, and R.~Krejcar.
\newblock Beta-integers as natural counting systems for quasicrystals.
\newblock {\em J. Phys. A}, 31(30):6449--6472, 1998.

\bibitem[Bla89]{Blanchard}
F.~Blanchard.
\newblock {$\beta$}-expansions and symbolic dynamics.
\newblock {\em Theoret. Comput. Sci.}, 65(2):131--141, 1989.

\bibitem[BS05]{BertheSiegel05}
V.~Berth{\'e} and A.~Siegel.
\newblock Tilings associated with beta-numeration and substitutions.
\newblock {\em Integers}, 5(3):A2, 46 pp. (electronic), 2005.

\bibitem[BS07]{BertheSiegel07}
V.~Berth{\'e} and A.~Siegel.
\newblock Purely periodic $\beta$-expansions in the {P}isot non-unit case.
\newblock To appear, 2007.

\bibitem[CF86]{CasselsFroehlich}
J.~W.~S. Cassels and A.~Fr{\"o}hlich, editors.
\newblock {\em Algebraic number theory}, London, 1986. Academic Press Inc.

\bibitem[CFS82]{CFS}
I.~P. Cornfeld, S.~V. Fomin, and Ya.~G. Sina{\u\i}.
\newblock {\em Ergodic theory}, volume 245 of {\em Grundlehren der
  Mathematischen Wissenschaften [Fundamental Principles of Mathematical
  Sciences]}.
\newblock Springer-Verlag, New York, 1982.
\newblock Translated from the Russian by A. B. Sosinski\u\i.

\bibitem[DKS96]{Solomyak}
K.~Dajani, C.~Kraaikamp, and B.~Solomyak.
\newblock The natural extension of the {$\beta$}-transformation.
\newblock {\em Acta Math. Hungar.}, 73(1-2):97--109, 1996.

\bibitem[Fro00]{FrougnyTemuco}
C.~Frougny.
\newblock Number representation and finite automata.
\newblock In {\em Topics in symbolic dynamics and applications (Temuco, 1997)},
  volume 279 of {\em London Math. Soc. Lecture Note Ser.}, pages 207--228.
  Cambridge Univ. Press, 2000.

\bibitem[FS92]{FrougnySolomyak92}
C.~Frougny and B.~Solomyak.
\newblock Finite beta-expansions.
\newblock {\em Ergodic Theory Dynamical Systems}, 12:45--82, 1992.

\bibitem[HI97]{HamaImahashi}
M.~Hama and T.~Imahashi.
\newblock Periodic {$\beta$}-expansions for certain classes of {P}isot numbers.
\newblock {\em Comment. Math. Univ. St. Paul.}, 46(2):103--116, 1997.

\bibitem[IR04]{ItoRao04}
S.~Ito and H.~Rao.
\newblock Purely periodic {$\beta$}-expansions with {P}isot unit base.
\newblock {\em Proc. Amer. Math. Soc.}, 133(4):953--964 (electronic), 2004.

\bibitem[Lot02]{FrougnyLothaire}
M.~Lothaire.
\newblock {\em Algebraic combinatorics on words}, volume~90 of {\em
  Encyclopedia of Mathematics and its Applications}.
\newblock Cambridge University Press, Cambridge, 2002.
\newblock (Chapter 7, written by C. Frougny).

\bibitem[Par60]{Parry60}
W.~Parry.
\newblock On the {$\beta $}-expansions of real numbers.
\newblock {\em Acta Math. Acad. Sci. Hungar.}, 11:401--416, 1960.

\bibitem[Pra99]{Praggastis}
B.~Praggastis.
\newblock Numeration systems and {M}arkov partitions from self-similar tilings.
\newblock {\em Trans. Amer. Math. Soc.}, 351(8):3315--3349, 1999.

\bibitem[QRY05]{QuRaoYang05}
Y.-H Qu, H.~Rao, and Y.-M. Yang.
\newblock Periods of {$\beta$}-expansions and linear recurrent sequences.
\newblock {\em Acta Arith.}, 120(1):27--37, 2005.

\bibitem[Rau88]{rau2}
G.~Rauzy.
\newblock Rotations sur les groupes, nombres alg\'ebriques, et substitutions.
\newblock In {\em S\'eminaire de Th\'eorie des Nombres, 1987--1988 (Talence,
  1987--1988)}. Univ. Bordeaux I, Talence., 1988.
\newblock Exp.\ No.\ 21, 12.

\bibitem[Roh61]{Rohlin}
V.~A. Rohlin.
\newblock Exact endomorphisms of a {L}ebesgue space.
\newblock {\em Izv. Akad. Nauk SSSR Ser. Mat.}, 25:499--530, 1961.

\bibitem[San02]{Sano02}
Y.~Sano.
\newblock On purely periodic beta-expansions of {P}isot numbers.
\newblock {\em Nagoya Math. J.}, 166:183--207, 2002.

\bibitem[Sch80]{Schmidt80}
K.~Schmidt.
\newblock On periodic expansions of {P}isot numbers and {S}alem numbers.
\newblock {\em Bull. London Math. Soc.}, 12(4):269--278, 1980.

\bibitem[Sie03]{Siegel03}
A.~Siegel.
\newblock Repr\'esentation des syst\`emes dynamiques substitutifs non
  unimodulaires.
\newblock {\em Ergodic Theory Dynam. Systems}, 23(4):1247--1273, 2003.

\bibitem[Sin06]{Sing06}
Bernd Sing.
\newblock Iterated function systems in mixed {E}uclidean and {$ p$}-adic
  spaces.
\newblock In {\em Complexus mundi}, pages 267--276. World Sci. Publ.,
  Hackensack, NJ, 2006.

\bibitem[ST07]{SiegelThuswaldner}
A.~Siegel and J.~M. Thuswaldner.
\newblock Topological properties of rauzy fractals.
\newblock preprint, 2007.

\bibitem[SW02]{SW02}
V.~F. Sirvent and Y.~Wang.
\newblock Self-affine tiling via substitution dynamical systems and {R}auzy
  fractals.
\newblock {\em Pacific J. Math.}, 206(2):465--485, 2002.

\bibitem[Thu89]{Thurston}
W.~P. Thurston.
\newblock Groups, tilings and finite state automata.
\newblock In {\em AMS Colloquium lectures}. AMS Colloquium lectures, 1989.

\bibitem[Thu06]{Thuswaldner06}
J{\"o}rg~M. Thuswaldner.
\newblock Unimodular {P}isot substitutions and their associated tiles.
\newblock {\em J. Th\'eor. Nombres Bordeaux}, 18(2):487--536, 2006.

\end{thebibliography}
